\documentclass[10pt,twoside]{article}
\usepackage{latexsym}
\usepackage{amsthm}
\usepackage{amsfonts}
\usepackage{amssymb}

\theoremstyle{plain}
  \newtheorem{theorem}{Theorem}[section]
  
  \newtheorem{proposition}{Proposition}[section]
  \newtheorem{lemma}{Lemma}[section]
  \newtheorem{corollary}{Corollary}[section]

\theoremstyle{remark}
  \newtheorem{remark}{Remark}[section]

\theoremstyle{definition}

\makeatletter

\@addtoreset{equation}{section}
\makeatother

\begin{document}

\title{
A priori estimates and weak solutions for the derivative nonlinear Schr\"{o}dinger equation on torus below $H^{1/2}$}
\author{Hideo Takaoka\thanks{This work was supported by JSPS KAKENHI Grant Number 10322794.}\\
Department of Mathematics, Hokkaido University\\
Sapporo, 060-0810, Japan\\
takaoka@math.sci.hokudai.ac.jp}

\markboth{Hideo Takaoka}{Derivative nonlinear Schr\"odinger equation on torus}
\pagestyle{myheadings}

\date{}
      
\maketitle

\begin{abstract}
We propose a priori estimates for a weak solution to the derivative nonlinear Schr\"odinger equation (DNLS) on torus with small $L^2$-norm datum in low regularity Sobolev spaces.
These estimates allow us to show the existence of solutions in $H^s(\mathbb{T})$ with some $s<1/2$ in a relatively weak sense.
Furthermore we make some remarks on the error estimates arising from the finite dimensional approximation solutions of DNLS using the Fourier-Lesbesgue type as auxiliary spaces, despite the fact that Nahmod, Oh, Rey-Bullet and Staffilani \cite{nors} have already seen them.

\end{abstract}

\noindent
{\it $2010$ Mathematics Subject Classification{\rm :}} 
Primary  35Q55;  Secondary 42B37.

\noindent
{\it Key words and phrases{\rm :}}
nonlinear Schr\"odinger equation; well-posedness.

\section{Introduction}\label{introduction-sec}
\indent

In this paper we consider the Cauchy problem for the derivative nonlinear Schr\"odinger equation (DNLS) on the torus:
\begin{eqnarray}\label{dnls-ori}
i\partial_tu+\partial_x^2u=-iu^2\partial_x\overline{u}-\frac{1}{2}|u|^4u+\mu[u]|u|^2u-\psi[u]u,\quad (t,x)\in [-T,T]\times\mathbb{T},
\end{eqnarray}
\begin{eqnarray}\label{data}
u(0,x)=u_0(x),\quad x\in\mathbb{T},
\end{eqnarray}
where
$$
\mu[u](t)=\frac{1}{2\pi}\int_0^{2\pi}|u(t,\theta)|^2\,d\theta,\quad
\psi[u](t)=\frac{1}{2\pi}\int_0^{2\pi}\left(2\mathrm{Im}(u\partial_x\overline{u})(t,\theta)-\frac{1}{2}|u|^4(t,\theta)\right)\,d\theta,
$$
$\mathbb{T}=\mathbb{R}/2\pi\mathbb{Z}$ is the torus and $u=u(t,x):[-T,T]\times\mathbb{T}\to \mathbb{C}$.
Our aim in this paper is to revisit the classical subject of a priori estimates of solutions to (\ref{dnls-ori})-(\ref{data}).
The equation (\ref{dnls-ori}) possesses important conserved quantities; conserved mass $M$, conserved energy $E$ , conserved momentum $P$, where
$$
M[u](t)=\int_0^{2\pi}|u(t,x)|^2\,dx,
$$
$$
E[u](t)=\int_0^{2\pi}\left(|\partial_xu(t,x)|^2-\frac{1}{2}\mathrm{Im}(|u(t,x)|^2u(t,x)\partial_x\overline{u(t,x)})-\frac{\mu}{2}|u(t,x)|^4\right)\,dx,
$$
$$
P[u](t)=\int_0^{2\pi}\left(\mathrm{Im}(\overline{u(t,x)}\partial_xu(t,x))+\frac12|u(t,x)|^4\right)\,dx.
$$
Formally solutions of (\ref{dnls-ori}) satisfy
$$
M[u](t)=M[u](0),\quad E[u](t)=E[u](0),\quad P[u](t)=P[u](0).
$$
Therefore $M[u],~E(u),~P[u]$ remain constant through time as well.

There are a handful of other form of the derivative nonlinear Schr\"odinger equations.
Supposing that $u$ is sufficiently smooth solution to (\ref{dnls-ori}), the $L^2$ norm remains constant through time (mass conservation)
$$
M[u](t)=\|u(t)\|_{L^2}^2=\|u_0\|_{L^2}^2.
$$
We impose by putting the constant $\mu=\|u_0\|_{L^2}^2/2\pi$ in (\ref{dnls-ori}), and define the transformation
\begin{eqnarray}\label{gauge-trans}
v(t,x)=e^{i{\cal G}[u](t,x+2\mu t)+\frac{i}{4}\mu^2t}u(t,x+2\mu t),
\end{eqnarray}
where
$$
{\cal G}[u](t,x)=\frac{1}{2\pi}\int_0^{2\pi}\int_{\theta}^x\left(|u(t,y)|^2-\mu\right)\,dyd\theta.
$$
This transformation is known as a gauge transformation.
By means of the gauge transformation, we can change the nonlinear terms in familiar type.
In \cite{herr}, if one rewrites the equation (\ref{dnls-ori}) as the transform (\ref{gauge-trans}), one obtains that $v=v(t,x)$ satisfies the simple one:
\begin{eqnarray}\label{dnls}
i\partial_tv+\partial_x^2v=i\partial_x(|v|^2v).
\end{eqnarray}
The two models (\ref{dnls-ori}) and (\ref{dnls}) are equivalent in some sense.
We point out that the nonlinearity $\partial_x(|v|^2v)$ is unfavorable.
The equation (\ref{dnls}) contains mixed derivative nonlinear terms $\partial_x(|v|^2v)=|v|^2v_x+2v^2\overline{v_x}$ relating only two terms $|v|^2v_x,~v^2\overline{v_x}$.
In dealing with the nonlinearity of the form $|v|^2v_x$, the standard energy method does not work and we encounter a difficulty of the derivative loss, see \cite{GH}.
In order to overcome this difficulty, we focus our attention on the equation (\ref{dnls-ori}), permitting more nonlinear terms than (\ref{dnls}).

Herr \cite{herr} proved the local well-posedness to (\ref{dnls-ori})-(\ref{data}) in $H^s$ for $s\ge 1/2$.
When the $L^2$ norm is small, one can combine the energy conservation law with local well-posedness theory to obtain the global well-posedness in $H^s$ for $s\ge 1$ (by the Gagliardo-Nirenberg interpolation inequality).
For global solutions below the energy threshold, the global well-posedness was obtained by Su Win \cite{su} to (\ref{dnls-ori})-(\ref{data}) in $H^s$ for $s>1/2$.
We remark that the index $s=1/2$ is well-posedness regularity threshold.
Indeed, the uniform continuity of the flow map, as a map from any ball of $H^s$ into $C([-T,T],\,H^s)$ at any time $T>0$, does not hold if $s<1/2$, see \cite{GH}.

Concerning the whole real line case without periodic boundary condition, the best local well-posedness was known in $H^s(\mathbb{R})$ for $s\ge 1/2$, see \cite{oz,ho1,ho2,ta}.
This result is sharp with respect to the lower threshold on $s$, which is essentially of the same kind as the one for the periodic boundary condition case.
Moreover there was also global well-posedness for data in $H^s(\mathbb{R})$ for $s\ge 1/2$, see \cite{ckstt1,my}.

In the present paper we will consider the existence of local in time solution in the case of periodic boundary condition for data below the threshold $s=1/2$.
More precisely, we prove the following theorem.
\begin{theorem}[existence and continuity of solutions]\label{thm-existence}
Let $4/9+a/9<s<\min\{1/2,3a/2\}$ and $a>8/25$.
There exists a positive constant $\varepsilon>0$ such that if $u_0\in H^s$ is small $L^2$ norm so that $\|u_0\|_{L^2}\le \varepsilon$, then there exist a positive time $T>0$ and a weak solution $u(t)$ to (\ref{dnls-ori})-(\ref{data}) on $t\in[-T,T]$ with
$$
u\in Y_T^a\cap C([-T,T],~H^s),
$$
$$
\sup_{-T\le t\le T}\|P_{\ge N}u(t)\|_{H^s}\le C\|P_{\ge CN}u_0\|_{H^s}+\frac{C}{N^{\varepsilon}},
$$
for any $N>0$, where $P_{\ge M}$ means the restriction operator with Fourier modes truncation $|\xi|\ge M$, and constants $C$ and $\varepsilon>0$ depend only on $\|u_0\|_{H^s}$.
Here the space $Y^a_T$ is defined in Section \ref{sec:2}.
\end{theorem}
\begin{remark}
The lower available value of $s$ achieved in Theorem \ref{thm-existence} is $s=12/25+4\varepsilon/9$ when $a=8/25+\varepsilon$.
\end{remark}
\begin{remark}
\begin{itemize}
\item[(i)]
In \cite{bialin}, Biagioni and Linares proved that the Cauchy problem associated to (\ref{dnls}) is ill-posed in $H^s$ for $s<1/2$ in the sense that the solution map fails to be uniformly continuous.
Furthermore, Gr\"unrock and Herr \cite{GH} mentioned that the failure of uniform continuity is shown in ${\cal F}L^{s,p}$ for $s<1/2$ and $r\in[1,\infty]$, where ${\cal F}L^{s,p}$ is called the Fourier-Lebesgue space,
$$
{\cal F}L^{s,p}=\{f\in {\cal D}'\mid \langle\xi\rangle^{s}\widehat{f}(\xi)\in \ell_{\xi}^p(\mathbb{Z})\},
$$
where 
$$
\|f\|_{{\cal F}L^{s,p}}=\|\langle\xi\rangle^{s}f(\xi)\|_{\ell_{\xi}^{p}},
$$
with $\ell_{\xi}^p$ denoting the standard $\ell^p$ sequence space.
\item[(ii)]
In \cite{GH}, Gr\"unrock and Herr proved that if $u_0\in {\cal F}L^{1/2,p}$ with $2<p<4$ then the local well-posedness holds.
One immediately sees that the solution $u(t)$ obtained in Theorem \ref{thm-existence} is unique under the additional assumption on the initial data $u_0\in {\cal F}L^{1/2,p}$ for some $2<p<4$\footnote{For more details, refer to the proof of Theorem \ref{approx} in Section \ref{sec:8} in this paper.}.
\end{itemize}
\end{remark}

The reminder of the paper contains the finite dimensional approximation result in a low regularity Sobolev spaces, which are essentially of the same kind as the one already obtained by Nahmod, Oh, Rey-Bullet and Staffilani in \cite{nors}.
It was shown in \cite{nors} that the dynamics of approximate that of the equation (\ref{dnls-ori}) in ${\cal F}L^{s,p}$ with  $s>1/2$ and $2<p<4$ along with the uniform probabilistic energy estimate for the approximating solutions had its origin in \cite{bour1} allows one to establish global well-posedness almost surely in ${\cal F}L^{s,p}$, where the key ingredient is the finite dimensional approximation lemma.
We revisit and deduce the strong approximation lemma in $H^s\cap{\cal F}L^{s_1,p}$ with $1/4<s<1/2<s_1$ and $2<p<4$.

Following \cite{nors}, consider the finite dimensional approximation of (DNLS):
\begin{eqnarray}
i\partial_tu^N+\partial_x^2u^N & = & -iP_{\le N}((u^N)^2\partial_x\overline{u^N})-\frac{1}{2}P_{\le N}(|u^N|^4u^N)\nonumber\\
& & +\mu[u^N] P_{\le N}(|u^N|^2u^N)-\psi[u^N]u^N,
\label{fdnls}
\end{eqnarray}
\begin{eqnarray}\label{fdata}
u^N(0,x)=P_{\le N}u_0(x),
\end{eqnarray}
where $P_{\le N}$ means the restriction operator with Fourier modes truncation $|\xi|\le N$. 
Comparing solutions of (\ref{dnls-ori})-(\ref{data}) and (\ref{fdnls})-(\ref{fdata}), we obtain a priori error estimates for the finite dimensional approximation.

Denote $\|u_0\|_{H^s\cap {\cal F}L^{s_1,p}}=\|u_0\|_{H^s}+\|u\|_{{\cal F}L^{s_1,p}}$.
The result is given by the following theorem.
\begin{theorem}[approximation lemma]\label{approx}
Let $1/4<s<1/2<s_1$ and $2<p<4$.
Let $N$ and $A$ be constants.
Assume that $u_0\in H^s\cap {\cal F}L^{s_1,p}$ be such that $\|u_0\|_{H^s\cap {\cal F}L^{s_1,\infty}}<A$, and the solution $u^N(t)$ of (\ref{fdnls}) with data (\ref{fdata}) satisfies the bound
\begin{eqnarray*}
\left\|u^N(t)\right\|_{H^s\cap {\cal F}L^{s_1,p}}\le A
\end{eqnarray*}
for all $t\in [-T,T]$ for some given $T>0$.
Then the Cauchy problem (\ref{dnls-ori}) with data (\ref{data}) is well-posed on $[-T,T]$ and there exists constants $C_j,~1\le j\le 3$ such that the solution $u(t)$ of (\ref{dnls-ori})-(\ref{data}) satisfies the following estimate
\begin{eqnarray}\label{app}
\left\|u(t)-u^N(t)\right\|_{H^{s'}\cap {\cal F}^{s_1',p}}
\le C_1 \exp [C_2(1+A)^{C_3}T]N^{\max\{s'-s,s_1'-s_1\}},
\end{eqnarray}
for all $t\in [-T,T],~1/4<s'<s$ and $1/2<s_1'<s_1$, provided the right-hand side of (\ref{app}) remains less than $1$.
\end{theorem}

As a byproduct of the a priori error estimates in Theorem \ref{approx}, we can prove almost global well-posedness for the initial data in the support of the canonical Gaussian measures on $H^s\cap {\cal F}L^{s_1,p}$ for each $1/4<s<1/2<s_1$ and $2<p<4$.
As it was explained before, this result was already known in \cite{nors}, where they proved that the local in time solutions can be extended to be global ones almost surely in ${\cal F}L^{s,p}$ for some $s>1/2$ and $2<p<4$.
Note that using Theorem \ref{approx} it is possible to give the a priori bound of $H^s$ norm of the solution to (\ref{dnls-ori}) as well as that of ${\cal F}L^{s_1,p}$ norm.
The proof of the almost sure global well-posedness in $H^s\cap {\cal F}^{s_1,p}$ is accomplished by using Theorem \ref{approx} based on the same argument as in \cite{nors}.
Hence we will only give a proof of Theorem \ref{approx} in this paper.

The outline of the paper is organized as follows.
In Section \ref{sec:2}, we give some additional notation that is used throughout the paper, and introduce the some dispersive properties of solutions of the linear Schr\"odinger equation.
In Section \ref{sec:3}, we divide the nonlinearity into "resonant" and "nonresonant" components.
In Sections \ref{sec:X^s} and \ref{mul-1}, we exploit several multilinear estimates.
In Section \ref{sec:6}, we derive the a priori estimates that are applied in Section \ref{sec:energy}.
In Section \ref{sec:energy}, we prove Theorem \ref{thm-existence}.
Finally, in Section \ref{sec:8}, we provide the proof of Theorem \ref{approx}.

\section{Notation and preliminary results}\label{sec:2}
\indent

In this section we define some notation that is used in this article, and present some preliminary results.

\subsection{Notation}
\indent

Let $\phi\in C^{\infty}_0(\mathbb{R})$ be a bump function adapted to $[-2,2]$ which equals to $1$ on $[-1,1]$.
Also define $\psi\in C^{\infty}(\mathbb{R})$ such that $\psi(\xi)=1-\phi(\xi)$.
Set $\phi_{\rho}(t)=\phi(\rho^{-1}t)$ and $\psi_{\rho}(\xi)=\psi(\rho^{-1}\xi)$ for $\rho>0$.

Let $\chi_T(t)$ be the characteristic function that is equal to $1$ on $|t|\le T$ and is equal to $0$ on $|t|>T$. 
For a set $A$, $1_A$ denotes the characteristic function of $A$. 

We prefer to use notation $\langle x\rangle =(1+|x|^2)^{1/2}$ for $x\in \mathbb{R}$. 
Write $\xi_{kl}$ for $\xi_k+\xi_l$.

The Fourier transform with respect to the space variable (discrete Fourier transform) is defined by
$$
{\cal F}_xf(\xi)=\frac{1}{\sqrt{2\pi}}\int_0^{2\pi}e^{-ix\xi}f(x)\,dx,\quad \xi\in\mathbb{Z},
$$
and with respect to the time variable by
$$
{\cal F}_tf(\tau)=\frac{1}{\sqrt{2\pi}}\int_{-\infty}^{\infty}e^{-it\tau}f(t)\,dt,\quad\tau\in\mathbb{R},
$$
and ${\cal F}={\cal F}_t{\cal F}_x$.
Particularly, the independent variable $t$ represents time, and thus $\tau$ is used for variable in time frequency space.
Therefore, $\xi$ will represent the Fourier transform variable with respect to space variable $x$.
We also use the same Fourier transform definitions $\widehat{u}(\xi)$ denote 
${\cal F}_xu(\xi)$, if the confusion does not arise from the above definition.  

For $1\le p,q\le \infty$, we use the mixed norm notation $\|f\|_{L^q_tL^p_x}$ with norm
$$
\|f\|_{L^q_tL^p_x}=\left(\int_{-\infty}^{\infty}\|f(t)\|_{L^p_x(\mathbb{T})}^q\,dt\right)^{1/q},
$$
with the obvious modification when $q=\infty$.
For $T>0$, we also use $\|f\|_{L^q_TL^p_x}$ to denote the norm
$$
\|f\|_{L^q_TL^p_x}=\left(\int_0^T\|f(t)\|_{L^p_x(\mathbb{T})}^q\,dt\right)^{1/q},
$$
with the obvious modification when $q=\infty$.

We use $c,~C$ to denote various constants, usually depending only on $s$.
We use $A\lesssim B$ to denote $A\le CB$ for some constant $C>0$.
Similarly, we write $A\sim B$ to mean $A\lesssim B$ and $B\lesssim A$.

For $N\in\mathbb{N}$, the operator $P_{\le N}$ denotes the restriction operator to the $N$ first Fourier modes, as is readily used.
The operators $P_{\ge N}$ and $P_{N}$ denote the restriction operators to $|\xi|\ge N$ and $|\xi|=N$ Fourier modes, respectively.

For $s,~b\in\mathbb{R}$ and $1\le p,~q\le\infty$, we define the ${\cal X}^{s,b}_{p,q}$ norm \cite{GH} by
$$
\|u\|_{{\cal X}^{s,b}_{p,q}}=\left\|\langle\xi\rangle^s\left\|\langle\tau+\xi^2\rangle^b{\cal F}u(\tau,\xi)\right\|_{L_{\tau}^q}\right\|_{\ell_{\xi}^p}.
$$
We will make use of two parameter spaces $X^{s,b}$ with norm \cite{bour0}
$$
\|u\|_{X^{s,b}}=\|u\|_{{\cal X}^{s,b}_{2,2}},
$$
and define the slightly stronger norm space $Y^s$ by
$$
\|u\|_{Y^s}=\|u\|_{{\cal X}^{s,1/2}_{2,2}}+\|u\|_{{\cal X}^{s,0}_{2,1}}.
$$
We also need the companion space $Z^s$ which is defined by the norm
$$
\|u\|_{Z^s}=\|u\|_{{\cal X}^{s,-1/2}_{2,2}}+\|u\|_{{\cal X}^{s,-1}_{2,1}}
$$
Also define the norm space ${\cal Y}^s$ by
$$
\|u\|_{{\cal Y}^{s,p}}=\|u\|_{{\cal X}^{s,1/2}_{p,2}}+\|u\|_{{\cal X}^{s,0}_{p,1}},
$$
and the relevant companion space ${\cal Z}^{s,p}$ by
$$
\|u\|_{{\cal Z}^{s,p}}=\|u\|_{{\cal X}^{s,-1/2}_{p,1}}+\|u\|_{{\cal X}^{s,-1}_{p,1}}.
$$
For $T>0$, we define the restriction norm spaces $X_T^{s,b}$
$$
X^{s,b}_T=\{u|_{-T\le t\le T}\mid u\in X^{s,b}\},
$$
with norm
$$
\|u\|_{X_T^{s,b}}=\inf\{\|U\|_{X^{s,b}}\mid U|_{-T\le t\le T}=u\}.
$$
Also define $Y^s_T,~Z^s_T,~{\cal Y}^{s,p}_T$ and ${\cal Z}^{s,p}_T$ in the same manner, respectively.
\begin{remark}\label{embed}
Using Riemann-Lebesgue lemma, we easily see that
$$
Y^s\hookrightarrow C(\mathbb{R};H^s),\quad {\cal Y}^{s,p}\hookrightarrow C(\mathbb{R};{\cal F}L^{s,p}). 
$$
\end{remark}
\begin{remark}[Lemma 3.2 in \cite{ckstt2}]\label{rem:dual}
We remark that there is a duality relationship between $Y^s$ and $Z^s$.
Indeed, one can verify that\footnote{Exactly the same proof in \cite{ckstt2} works for the Schr\"odinger equation, while the KdV equation is considered in \cite{ckstt2}.}
\begin{eqnarray}\label{duality-Y^s}
\left|\int_{\mathbb{R}\times\mathbb{T}}\chi_{T}(t)u(t,x)
v(t,x)\,dtdx\right|\lesssim \|u\|_{Y^s}\|v\|_{Z^{-s}},
\end{eqnarray}
for all $s\in\mathbb{R}$ and $T>0$.
In particular, if $0<T_1<T_2$, then
\begin{eqnarray}\label{t-inc}
\|u\|_{X^{s,b}_{T_1}}\le \|u\|_{X^{s,b}_{T_2}}.
\end{eqnarray}
\end{remark}
For complex-valued $n$ functions $f_1,~f_2,\ldots,~f_n$ defined on the set $\mathbb{Z}$ of integers, we write
the discrete convolution (convolution sum) $[f_1*f_2*\ldots *f_n](\xi)$ as
$$
[f_1*f_2*\ldots *f_n](\xi)=\sum_{*}\prod_{j=1}^nf_1(\xi_j),
$$
where $\sum_*$ denotes a summation over the set where $\xi_1+\xi_2+\ldots +\xi_n=\xi$.
Also write
$$
[g_1*g_2*\ldots*g_n](\tau)=\int_*\prod_{j=1}^ng_j(\tau_j),
$$
where $\int_*$ denotes an integration over the set where $\tau_1+\tau_2+\ldots+\tau_n=\tau$.

It is convenient to introduce some useful notation for multilinear expressions.
If $k\ge 2$ is an even integer, we define the hyperplane
$$
\Gamma_n=\{(\xi_1,\ldots,\xi_n)\in\mathbb{R}^n\mid \xi_1+\ldots+\xi_n=0\}.
$$
For any function $m(\xi_1,\ldots,\xi_n)$ on $\Gamma_n$, we define the $n$-multilinear form
$$
\sum_{\Gamma_n}m(\xi_1,\ldots,\xi_n)\prod_{j=1}^nf(\xi_j)=\sum_{(\xi_1,\ldots,\xi_n)\in \Gamma_n}m(\xi_1,\ldots,\xi_n)\prod_{j=1}^nf(\xi_j).
$$
Also define
$$
\int_{\Gamma_n}m(\tau_1,\ldots,\tau_n)\prod_{j=1}^ng(\tau_j)=\int_{\tau_1+\ldots+\tau_n=0}\left(\prod_{j=1}^ng(\tau_j)\right)\,d\tau_1\ldots d\tau_{n-1}.
$$

\subsection{Dispersive estimates}
\indent

In this subsection, we list a series of estimates for solutions of linear problem and the inhomogeneous problem associated to the equation (\ref{dnls-ori}).
\begin{lemma}\label{linear}
For all $s\in\mathbb{R}$,
\begin{eqnarray}\label{l1}
\|\phi(t)e^{it\partial_x^2}u_0\|_{Y^s}\lesssim \|u_0\|_{H^s},
\end{eqnarray}
\begin{eqnarray}\label{l2}
\left\|\phi(t)\int_0^te^{i(t-t')\partial_x^2}f(t')\,dt'\right\|_{Y^s}\lesssim\|f\|_{Z^s}.
\end{eqnarray}
\end{lemma}
\noindent
{\it Proof.}
For (\ref{l1}) and (\ref{l2}), see \cite[Lemmas 7.1 and 7.2]{CKSTT}.

\begin{lemma}\label{lem:L^pL^q}
\begin{itemize}
\item[(i)]
If $2\le p<\infty,~b\ge 1/2-1/p$, we have
\begin{eqnarray}\label{l3}
\|u\|_{L_t^pH_x^s}\lesssim \|u\|_{X^{s,b}}.
\end{eqnarray}
\item[(ii)]
If $2\le p,q<\infty,~b\ge 1/2-1/p,~s\ge 1/2-1/q$, we have 
\begin{eqnarray}\label{l4}
\|u\|_{L_t^pL_x^q}\lesssim \|u\|_{X^{s,b}}.
\end{eqnarray}
\item[(iii)]
If $1<p\le 2,~b\le 1/2-1/p$, we have
\begin{eqnarray}\label{l5}
\|u\|_{X^{s,b}}\lesssim \|u\|_{L_t^pH_x^s}.
\end{eqnarray}
\item[(iv)]
If $-b',b>3/8$, we have 
\begin{eqnarray}\label{l6}
\|u\|_{L_{t,x}^4}\lesssim \|u\|_{X^{0,b}},
\end{eqnarray}
and
\begin{eqnarray}\label{l7}
\|u\|_{X^{0,b'}}\lesssim \|u\|_{L_{t,x}^{4/3}}.
\end{eqnarray}
\end{itemize}
\end{lemma}
\noindent
{\it Proof.}
See \cite{bour0} and \cite{herr}.
\begin{remark}
Interpolating between (\ref{l4}) with $p=8,~q=2$ and (\ref{l6}), we have that
\begin{eqnarray}\label{ip1}
\|u\|_{L_t^{\frac{8}{2-\varepsilon}}L_x^{\frac{4}{1+\varepsilon}}}\lesssim \|u\|_{X^{0,b}},
\end{eqnarray}
for $0<\varepsilon<1$ and $b>3/8$.
Also by (\ref{l4}) we have 
\begin{eqnarray}\label{ip2}
\|u\|_{L_t^{\frac{4}{\varepsilon}}L_x^{\frac{4}{1-2\varepsilon}}}\lesssim \|u\|_{X^{\frac{1+2\varepsilon}{4},b}},
\end{eqnarray}
for $0<\varepsilon<1/2$ and $b>1/2-\varepsilon/4$.
\end{remark}

\begin{lemma}\label{l-gain}
Let $s\in\mathbb{R}$ and $0<T<1$.
\begin{itemize}
\item[(i)]
For $0\le b_1<b_2<1/2$ or $-1/2<b_1<b_2\le 0$, there exists $c>0$ such that
\begin{eqnarray}\label{l8}
\|\phi_Tf\|_{X^{s,b_1}}\le cT^{b_2-b_1}\|f\|_{X^{s,b_2}}.
\end{eqnarray}
\item[(ii)]
For any $\delta>0$ there exists $c>0$ such that
\begin{eqnarray}\label{l9}
\|\phi_Tf\|_{X^{s,1/2}}\le cT^{-\delta}\|f\|_{X^{s,1/2}},
\end{eqnarray}
and
\begin{eqnarray}\label{l10}
\|\phi_Tf\|_{Y^{s}}\le cT^{-\delta}\|f\|_{Y^s}.
\end{eqnarray}
\item[(iii)]
For $0<b<1/2$, there exists $c>0$ such that
\begin{eqnarray}\label{l11}
\|\chi_Tf\|_{X^{s,b}}\le c\|f\|_{X^{s,b}}.
\end{eqnarray}
\end{itemize}
\end{lemma}
\noindent
{\it Proof.}
See \cite{herr} for the proof of (i) and (ii).
The proof of (iii) follows from the Leibniz rule for fractional derivative, $\chi_T(t)\in H_t^{b}$ and $\|\chi_T\|_{H^b}\le c$.
\qed

\begin{remark}\label{replace}
All of the estimates (\ref{l1})-(\ref{l2}) in Lemma \ref{linear} and (\ref{l8})-(\ref{l11}) in Lemma \ref{l-gain} still hold with ${\cal Y}^{s,p},~{\cal Z}^{s,p},~{\cal F}L^{s,p}$ replacing $Y^s,~Z^s,~H^s$, respectively.
For the proof, see \cite[Lemma 7.1]{GH}.
\end{remark}

In \cite{GH}, the following trilinear estimate was proven.
\begin{lemma}\cite[Lemma 5.1]{GH}
For $1/3<b<1/2$ and $s>3(1/2-b)$, the estimate
\begin{eqnarray}\label{GH-trilinear}
\left\|\prod_{J=1}^3u_j\right\|_{L_{t,x}^2}\lesssim \|u_1\|_{X^{s,b}}\|u_2\|_{X^{s,b}}\|u_3\|_{X^{0,b}}
\end{eqnarray}
holds true.
\end{lemma}

\section{A resonant decomposition}\label{sec:3}
\indent

In this section we discuss the structural nonlinear properties of the equation (\ref{dnls-ori}).
Defining
$$
N[u]=-iu^2\partial_x\overline{u}-\frac{1}{2}|u|^4u+\mu|u|^2u-\psi[u]u,
$$
\begin{eqnarray*}
N_1[u]=-iu^2\partial_x\overline{u}-\frac{1}{2\pi}\left(\int_0^{2\pi}2\mathrm{Im}(u\partial_x\overline{u})(t,\theta)\,d\theta\right) u,
\end{eqnarray*}
and
$$
N_2[u]=-\frac{1}{2}|u|^4u+\mu|u|^2u+\frac{1}{4\pi}\left(\int_0^{2\pi}|u(t,\theta)|^4\,d\theta\right)u,
$$
where $\mu=\|u_0\|_{L^2}^2/2\pi$ is a constant, we have that the nonlinear term $N[u]$ in equation (\ref{dnls-ori}) can be decomposed of an effective cubic nonlinear term with derivative $N_1[u]$, plus other terms $N_2[u]$ without spatial derivatives.
In this section we reformulate the cubic derivative nonlinear term $N_1[u]$ with a resonant decomposition.

The derivative cubic nonlinear terms $N_1[u]$ are roughly classified into a nonlinear regimes of the non-resonance interaction modes and the resonance interaction modes.
Firstly, we identify these interaction modes.
The reason behind this classification is that the resonance interaction mode can be easier to handle by taking energy estimates in subsection \ref{H_x^s-estimates}.

We adapt the spatial Fourier transform to $N_1[u]$, so that
\begin{eqnarray*}
\widehat{N_1[u]}(t,\xi)=& & \frac{1}{2\pi}\sum_{*}\widehat{u}(t,\xi_1)\xi_2\widehat{\overline{u}}(t,\xi_2)\widehat{u}(t,\xi_3)\\
& & -\frac{1}{\pi}\left(\sum_{\xi_1+\xi_2=0}\xi_2\widehat{u}(t,\xi_1)\widehat{\overline{u}}(t,\xi_2)\right)\widehat{u}(t,\xi).
\end{eqnarray*}
Let us consider the algebraic identity: for $\xi=\xi_1+\xi_2+\xi_3$
$$
\xi_1^2-\xi_2^2+\xi_3^2-\xi^2=-2(\xi_1-\xi)(\xi_3-\xi).
$$
Using this identity we distinguish summation over all indices $\xi,~\xi_j~(1\le j\le 3)$;
\begin{itemize}
\item[(i)]
$(\xi_1-\xi)(\xi_3-\xi)\ne 0$,
\item[(ii)]
$\xi_3=\xi,~\xi_1+\xi_2=0$,
\item[(iii)]
$\xi_1=\xi,~\xi_2+\xi_3=0$.
\end{itemize}
The cases (ii) and (iii) are not complementary to each other.
The case for redundancy between (ii) and (iii) is $\xi_1=\xi_3=-\xi_2=\xi$.
Due to the fact that $(\xi_1-\xi)(\xi_3-\xi)=0$ in (ii) or (iii), we have
$$
\sum_{\scriptstyle \xi_1+\xi_2+\xi_3=\xi \atop{\scriptstyle (\xi_1-\xi)(\xi_3-\xi)=0}}=\sum_{\scriptstyle \xi_3=\xi \atop{\scriptstyle \xi_1+\xi_2=0}}+\sum_{\scriptstyle \xi_1=\xi \atop{\scriptstyle \xi_2+\xi_3=0}}-\sum_{\scriptstyle \xi_1=\xi_3=-\xi_2=\xi \atop{\scriptstyle \xi_1+\xi_2=\xi_2+\xi_3=0}}.
$$
Because of this, we have
\begin{eqnarray*}
\widehat{N_1[u]}(t,\xi)& = & \frac{1}{2\pi}\sum_{\scriptstyle * \atop{\scriptstyle (\xi_1-\xi)(\xi_3-\xi)\ne 0}}\widehat{u}(t,\xi_1)\xi_2\widehat{\overline{u}}(t,\xi_2)\widehat{u}(t,\xi_3)+\frac{1}{2\pi}\xi|\widehat{u}(t,\xi)|^2\widehat{u}(t,\xi)\\
& = & \widehat{N_{11}[u]}(t,\xi)+\widehat{N_{12}[u]}(t,\xi),
\end{eqnarray*}
where $N_{11}[u]=N_{11}(u,u,u),~N_{12}[u]=N_{12}(u,u,u)$,
\begin{eqnarray}\label{T_1}
\widehat{N_{11}(u_1,u_2,u_3)}(t,\xi)=\frac{1}{2\pi}\sum_{\scriptstyle * \atop{\scriptstyle (\xi_1-\xi)(\xi_3-\xi)\ne 0}}\widehat{u_1}(t,\xi_1)\xi_2\widehat{\overline{u_2}}(t,\xi_2)\widehat{u_3}(t,\xi_3),
\end{eqnarray}
and
\begin{eqnarray}\label{T_2}
\widehat{N_{12}(u_1,u_2,u_3)}(t,\xi)=\frac{1}{2\pi}\xi\widehat{u_1}(t,\xi)\widehat{\overline{u_2}}(t,-\xi)\widehat{u_3}(t,\xi).
\end{eqnarray}
Since the nonlinear resonance forced by $(\xi_1-\xi)(\xi_3-\xi)= 0$ is the occurrence resonance of resonance in a nonlinearity $N_1[u]$,  we say $N_{11}[u]$ and $N_{12}[u]$ as the non-resonance and resonance terms, respectively.
The resonance term $N_{12}[u]$ corresponds to forced oscillations that may oscillate with greater amplitude than at $N_{11}[u]$.

For other terms in $N_1[u]$, define
\begin{eqnarray}\label{N_{21}}
N_{21}[u]=\mu|u|^2u,
\end{eqnarray}
$$
N_{22}[u]=-\frac{1}{2}|u|^4u
+\frac{1}{8\pi}\left(\int_0^{2\pi}|u(t,\theta)|^4\,d\theta\right) u,
$$
so that
$$
N_2[u]=N_{21}[u]+N_{22}[u].
$$
In conclusion, we show that the nonlinear term $N[u]$ of equation can be expanded as follows:
\begin{eqnarray}\label{dnls-ori3}
N[u] = N_1[u]+N_2[u]=\sum_{k,l=1}^2 N_{kl}[u].
\end{eqnarray}
\begin{remark}
We shall need $s\ge 1/2$, if we are to control terms $N_{11}[u]$ and $N_{12}[u]$ by the Picard iteration scheme on the integral equation associated to (\ref{dnls-ori})-(\ref{data}).
Indeed, in \cite{herr} it is shown that the Cauchy problem (\ref{dnls})-(\ref{data}) is analytically locally well-posed in $H^s$ for $s\ge 1/2$.
Moreover in \cite{bialin}, the Cauchy problem (\ref{dnls})-(\ref{data}) is shown to be locally ill-posed in $H^s$ for $s<1/2$.
The key estimate in which local well-posedness in $H^s$ for $s\ge 1/2$ is the trilinear $X^{s,b}$ estimate
$$
\|u_1\partial_x\overline{u_2}u_3\|_{X^{s,-1/2}}\le \|u_1\|_{X^{s,1/2}}\|u_2\|_{X^{s,1/2}}\|u_3\|_{X^{s,1/2}}
$$
for all functions $u_j~(1\le j\le 3)$, where $s\ge 1/2$.
Using this and the standard computation \cite{bour1}, we obtain the local existence theory for $s\ge 1/2$.
To get down to $s\ge 1/2$, we prove local a priori estimates for energy-based methods.
\end{remark}

\section{Multilinear estimates I}\label{sec:X^s}
\noindent

In this section we illustrate several multilinear estimates.
\subsection{Trilinear estimates}\label{sec:4}
\noindent

We take the advantage of the identity
\begin{eqnarray}\label{id}
\sum_{j=1}^4\left(\tau_j+(-1)^{j-1}\xi_j^2\right)=2\xi_{12}\xi_{14}
\end{eqnarray}
which holds whenever $\sum_{j=1}^4\tau_j=\sum_{j=1}^4\xi_j=0$.

\begin{lemma}\label{thm:tri-linear}
Let $4/9+a/9<s<1/2$ and $a>1/4$.
Then there exists $1/3<b<1/2$ such that
\begin{eqnarray*}
& & \|N_{11}(u_1,u_2,u_3)\|_{X^{a,-1/2}}\\
& \lesssim & \sum_{k=1}^3\|u_k\|_{X^{a,1/2}}\prod_{j=1,\ne k}^3\|u_j\|_{X^{a,b}}+ \sum_{k=1}^3\|u_k\|_{X^{a,1/2}}\prod_{j=1,\ne k}^3\|u_j\|_{L^8_tH_x^s}\\
& & 
+\sum_{k=1}^3\|u_k\|_{L_t^{2}H_x^s}\prod_{j=1,\ne k}^3\left(\|u_j\|_{L_t^{\infty}H_x^s}+\|u_j\|_{L^8_tH_x^s}\right).
\end{eqnarray*}
\end{lemma}
\noindent
{\it Proof.}
We require the following estimates:
\begin{eqnarray}
& & \left\|\langle \xi\rangle^a\frac{{\cal F}N_{11}(u_1,u_2,u_3)(\tau,\xi)}{\langle\tau+\xi^2\rangle^{1/2}}\right\|_{L_{\tau}^2\ell_{\xi}^2}\label{N-11-1}\\
& \lesssim &\sum_{k=1}^3\|u_k\|_{X^{a,1/2}}\prod_{j=1,\ne k}^3\|u_j\|_{X^{a,b}}+ \sum_{k=1}^3\|u_k\|_{X^{a,1/2}}\prod_{j=1,\ne k}^3\|u_j\|_{L^8_tH_x^s}\nonumber\\
& & 
+\sum_{k=1}^3\|u_k\|_{L_t^{2}H_x^s}\prod_{j=1,\ne k}^3\left(\|u_j\|_{L_t^{\infty}H_x^s}+\|u_j\|_{L^8_tH_x^s}\right).\nonumber
\end{eqnarray}
By (\ref{T_1}), we see that
$$
{\cal F}N_{11}(u_1,u_2,u_3)(\tau,\xi)=c\sum_{\scriptstyle * \atop{\scriptstyle (\xi_1-\xi)(\xi_3-\xi)\ne 0}}\int_*{\cal F}u_1(\tau_1,\xi_1){\cal F}\overline{u_2}(\tau_2,\xi_2){\cal F}u_3(\tau_3,\xi_3).
$$
Use the dyadic partition 
$$
N_j\sim\langle\xi_j\rangle,~K_j\sim\langle\tau_j+(-1)^{j-1}\xi_j^2\rangle\quad\mbox{for $1\le j\le 3$},$$
$$
N_4\sim\langle\xi\rangle,~K_4\sim\langle\tau+\xi^2\rangle,
$$
and
$$
N_{12}\sim|\xi_{12}|,~N_{14}\sim|\xi_{14}|.
$$
Using Littlewoods-Paley decomposition for $u_j$, we separate the integral and sum of the areas into the following cases:
\begin{itemize}
\item[$(A_1)$]
$N_2\ll\min\{N_1,N_3\}$ and $\max\{K_1,K_3\}= \max\{K_1,K_2,K_3,K_4\}$,
\item[$(A_2)$]
$\min\{N_1,N_3\}\lesssim N_2\lesssim \max\{N_1,N_3\}$ and $N_4\ll N_2$,
\item[$(A_3)$]
$N_2\gg \max\{N_1,N_3\}$ and $\max\{K_2,K_4\}=\max\{K_1,K_2,K_3,K_4\}$,
\item[$(A_4)$]
$N_{2}\ll\min\{N_1,N_3\}$ and $\max\{K_2,K_4\}= \max\{K_1,K_2,K_3,K_4\}$,
\item[$(A_5)$]
$\min\{N_1,N_3\}\lesssim N_2\lesssim\max\{N_1,N_3\}$ and $N_{12}\gtrsim N_4\gtrsim N_2$,
\item[$(A_6)$]
$N_2\gg \max\{N_1,N_3\}$ and $\max\{K_1,K_3\}=\max\{K_1,K_2,K_3,K_4\}$,
\item[$(A_7)$]
$\min\{N_1,N_3\}\lesssim N_2\lesssim\max\{N_1,N_3\}$ and $N_2\lesssim N_{12}\ll N_4$,
\item[$(A_8)$]
$\min\{N_1,N_3\}\lesssim  N_2\lesssim\max\{N_1,N_3\}$ and $\max\{N_2,N_4\}\gg \min\{N_2,N_4\}\gg N_{12}$,
\item[$(A_9)$]
$\min\{N_1,N_3\}\lesssim N_2\lesssim\max\{N_1,N_3\}$ and $\max\{N_2,N_4\}\sim \min\{N_2,N_4\}\gg N_{12}$.
\end{itemize}
In cases $(A_j),~1\le j\le 8$, we estimate the contributions of these cases to the left-hand side of (\ref{N-11-1}) by
$$
c\sum_{k=1}^3\|u_k\|_{X^{a,1/2}}\prod_{j=1,\ne k}^3\|u_j\|_{X^{a,b}}.
$$
On the other hand, in case $(A_9)$, we estimate the contribution of this case to the left-hand side of (\ref{N-11-1}) by
$$
c\sum_{k=1}^3\|u_k\|_{X^{a,1/2}}\prod_{j=1,\ne k}^3\|u_j\|_{L^8_tH_x^s}
+c\sum_{k=1}^3\|u_k\|_{L_t^{2}H_x^s}\prod_{j=1,\ne k}^3(\|u_j\|_{L_t^{\infty}H_x^s}+\|u_j\|_{L^8_tH_x^s}).
$$
We will postpone the proof of case $(A_9)$ in Lemma \ref{special case}, and consider the cases that from $(A_1)$ through $(A_8)$ here. 

In cases $(A_j),~1\le j\le 8$, it is convenient to use the change of variables $\xi_4=-\xi$ and $\tau_4=-\tau$.
Using duality, in these cases, it suffices to prove that
\begin{eqnarray}
& &\left| \int_{\Gamma_4}\sum_{\Gamma_4}1_{A_j,\xi_{14}\xi_{34}\ne 0}(\overline{\tau},\overline{\xi}){\cal F}{u_1}(\tau_1,\xi_1){\cal F}{\overline{\partial_xu_{2}}}(\tau_2,\xi_2) {\cal F}{u_3}(\tau_3,\xi_3){\cal F}{\overline{w}}(\tau_4,\xi_4)\right|\label{tri_X}\\
 &\lesssim  & \|w\|_{X^{-a,1/2}}\sum_{k=1}^3\|u_k\|_{X^{a,1/2}}\prod_{j=1,\ne k}^3\|u_j\|_{X^{a,b}},\nonumber
\end{eqnarray}
for $1\le j\le 8$, where
$$
\overline{\tau}=(\tau_1,\tau_2,\tau_3,\tau_4),~\overline{\xi}=(\xi_1,\xi_2,\xi_3,\xi_4).
$$
Moreover putting
$$
{\cal F}{v_j}(\tau,\xi)=\langle\xi\rangle^a{\cal F}{u_j}(\tau,\xi),
$$
for $1\le j\le 3$, and
$$
{\cal F}{v_4}(\tau,\xi)=\langle\xi\rangle^{-a}{\cal F}{w}(\tau,\xi),
$$
we can write (\ref{tri_X}) as the following equivalent form
\begin{eqnarray}
& &\left|\int_{\Gamma_4}\sum_{\Gamma_4} M(\overline{\xi},\overline{\tau}){\cal F}{v_1}(\tau_1,\xi_1){\cal F}{\overline{v_2}}(\tau_2,\xi_2) {\cal F}{v_3}(\tau_3,\xi_3){\cal F}{\overline{v_4}}(\tau_4,\xi_4)\right|\label{i2}\\
& \lesssim & \prod_{j=1}^4\|v_j\|_{X^{0,1/2}},\nonumber
\end{eqnarray}
where the multiplier $M$ is of the form
$$
M(\overline{\tau},\overline{\xi})=1_{A_j,\xi_{14}\xi_{34}\ne 0}(\overline{\tau},\overline{\xi})\frac{\langle\xi_2\rangle^{1-a}\langle\xi_4\rangle^a}{\langle\xi_1\rangle^a\langle\xi_3\rangle^a}.
$$
The definition of the norm $X^{s,b}$ allows one to keep $\widehat{u}$ and $\widehat{w}$ by positive functions.

\underline{Case $(A_1)$.}
In this case, we see that
\begin{eqnarray}\label{i-1}
|\xi_{12}|\sim N_1,~|\xi_{14}|\sim N_3,~N_4\lesssim\max\{N_1,N_3\}.
\end{eqnarray}
By symmetry we may assume $K_1\ge K_3$.
By (\ref{id}), we see that
$$
\langle\tau_1+\xi_1^2\rangle \gtrsim |\xi_{12}||\xi_{14}|\sim |\xi_1||\xi_3|.
$$
Using $N_2\ll\min\{N_1,N_3\}$ and $N_4\lesssim \max\{N_1,N_3\}$, we have the bound
$$
M(\overline{\tau},\overline{\xi})\lesssim
\frac{\langle\tau_1+\xi_1^2\rangle^{1/2}}{\langle\xi_3\rangle^{2a-1/2}\langle\xi_4\rangle^{2a-1/2}}.
$$
We group $v_2,v_3,v_4$ together and apply (\ref{GH-trilinear}) to control the contribution of $(A_1)$ to the left-hand side of (\ref{i2}) by
$$
c\|v_1\|_{X^{0,1/2}}\prod_{j=2}^4\|v_j\|_{X^{0,b}},
$$
for $a>1/4$ and $b>1/2-(2a-1/2)/3$, which has the desired estimate.

\underline{Case $(A_2)$.}
We show that the estimate corresponding to (\ref{i2}) with replacing the integral and sum of areas that contributions by case $(A_2)$.
Symmetry properties permits us to assume $N_3\lesssim N_2\lesssim N_1$.
Since $N_4\ll N_2$, we have $N_1\sim N_2$ and $ |\xi_{14}|\sim N_1$.
Therefore we have
$$
M(\overline{\tau},\overline{\xi})\lesssim N_1^{\frac12-2a}\max\{K_1,K_2,K_3,K_4\}^{1/2}.
$$
We compress this bounds into the discussion presented in case $(A_1)$, which shows that the contribution of the left-hand side of (\ref{i2}) is bounded by
$$
c\sum_{k=1}^3\|v_k\|_{X^{0,1/2}}\prod_{j=1,\ne k}^3\|v_j\|_{X^{0,b}},
$$
provided $a>1/4$ and $b>1/2-(2a-1/2)/3$.

\underline{Case $(A_3)$.}
In this case, we have $N_4\sim N_2,~|\xi_{12}|\sim N_2$ and $|\xi_{14}|\sim N_2$ which implies $\max\{K_2,K_4\}\gtrsim N_2^2$.
By symmetry we may assume $K_4\ge K_2$.
Then
$$
M(\overline{\tau},\overline{\xi})\lesssim \frac{K_4^{1/2}}{N_1^aN_3^a}.
$$
Therefore, we are recast the estimate (\ref{i2}) by grouping $v_1,v_2,v_3$ with the $L_{t,x}^2$ estimate given by (\ref{GH-trilinear}) as
$$
c\|v_4\|_{X^{0,1/2}}\prod_{j=1}^3\|v_j\|_{X^{0,b}},
$$
provided $a>0$ and $b>1/2-a/3$.

\underline{Case $(A_4)$.}
This condition implies $N_4\lesssim \max\{N_1,N_3\}$.
The argument analogous to proof of case $(A_3)$ shows that
$$
M(\overline{\tau},\overline{\xi})\lesssim
\frac{\max\{K_2,K_4\}^{1/2}}{\min\{N_1,N_3\}^{2a-1/2}\max\{N_1,N_3\}^{1/2}}.
$$
By symmetry we may assume $K_4\ge K_2$.
Therefore, the contribution of left-hand side of (\ref{i2}) to this case is estimated, via (\ref{GH-trilinear}), by
$$
c \|v_4\|_{X^{0,1/2}}\prod_{j=1}^3\|v_j\|_{X^{0,b}},
$$
provided $a>0$ and $b>1/2-(2a-1/2)/3$.

\underline{Case $(A_5)$.}
In the subregion where $N_1\ge N_3$, the convolution constrain $\sum_{j=1}^4\xi_j=0$ implies $\max\{N_3,N_{14}\}\gtrsim N_2$.
Then
\begin{eqnarray*}
M(\overline{\tau},\overline{\xi})& \lesssim & 
\frac{N_2^{\frac12-a+2\varepsilon}}{N_1^{\varepsilon}N_3^aN_4^{\varepsilon}N_{14}^{\frac12}}\max\{K_1,K_2,K_3,K_4\}^{1/2}\\
& \lesssim &\frac{\max\{K_1,K_2,K_3,K_4\}^{1/2}}{\max\{N_1,N_2,N_3,N_4\}^{2\varepsilon}N_2^{2a-\frac12-4\varepsilon}},
\end{eqnarray*}
for $0<\varepsilon\le 1/2-a$ sufficiently small.

In the subregion when $N_1<N_3$, it follows that $N_2^{1-a}N_4^a\le N_2^aN_4^{1-a}$.
Then the same argument as above shows that
\begin{eqnarray*}
M(\overline{\tau},\overline{\xi})& \lesssim & 
\frac{N_4^{\frac12-a+2\varepsilon}}{N_3^{\varepsilon}N_1^aN_2^{\varepsilon}N_{14}^{\frac12}}\max\{K_1,K_2,K_3,K_4\}^{1/2}\\
& \lesssim &\frac{\max\{K_1,K_2,K_3,K_4\}^{1/2}}{\max\{N_1,N_2,N_3,N_4\}^{2\varepsilon}N_4^{2a-\frac12-4\varepsilon}},
\end{eqnarray*}
for $0<\varepsilon\le 1/2-a$ sufficiently small.
As a consequence, the contribution of the left-hand side of (\ref{i2}) to this region, via using (\ref{GH-trilinear}), has the desired estimate, where $a>1/4$ and $b>1/2-2\varepsilon/3$.

\underline{Case $(A_6)$.}
In this case, we have $N_2\sim N_4\gg \max\{N_1,N_3\}$, which implies $N_{12}\sim N_{14}\sim N_2$.
Then
$$
M(\overline{\tau},\overline{\xi})\lesssim \frac{\max\{K_1,K_3\}^{1/2}}{N_1^aN_3^a}.
$$
By symmetry we may suppose $K_1\ge K_3$.
We group $v_2\overline{v_3}v_4$ in $L_t^2L_x^{\frac43}$ and estimate the contribution of this case to the left-hand side of (\ref{i2}) by
$$
c\|\langle D_x\rangle^{-a}w_1\|_{L_{t}^2L_x^{4}}\|\overline{v_2}(\langle D_x^{-a} \rangle v_3)\overline{v_4}\|_{L_t^2L_x^{\frac43}},
$$
where
$$
{\cal F}w_1(\tau,\xi)=\langle \tau+\xi^2\rangle^{1/2}{\cal F}v_1(\tau,\xi).
$$
Using the Sobolev inequality, we have that
$$
\|\langle D_x\rangle^{-a}w_1\|_{L_t^2L_x^{4}}\lesssim \|w_1\|_{L_{t,x}^2}=\|v_1\|_{X^{0,1/2}}
$$
for $a>1/4$.
On the other hand, by (\ref{ip1}) and (\ref{ip2}), it follows that
\begin{eqnarray*}
\|\overline{v_2}(\langle D_x^{-a} \rangle v_3)\overline{v_4}\|_{L_t^2L_x^{\frac43}}& \lesssim &
\|v_2\|_{L_{t}^{\frac{8}{2-\varepsilon}}L_x^{\frac{4}{1+\varepsilon}}}\|\langle D_x\rangle^{-a}v_3\|_{L_{t}^{\frac{4}{\varepsilon}}L_x^{\frac{4}{1-2\varepsilon}}}\|v_4\|_{L_{t}^{\frac{8}{2-\varepsilon}}L_x^{\frac{4}{1+\varepsilon}}}\\
& \lesssim & 
\prod_{j=2}^4\|v_j\|_{X^{0,b}},
\end{eqnarray*}
for $a>1/4+\varepsilon/2,~b>1/2-\varepsilon/4$ and $0<\varepsilon<1/2$.
Then the desired estimate follows in this case.

\underline{Case $(A_7)$.}
In this case region, we observe that $N_4\sim N_{14}$ which implies
$$
M(\overline{\tau},\overline{\xi})\lesssim \frac{\max\{K_1,K_2,K_3,K_4\}^{1/2}}{N_1^aN_2^{a-\varepsilon}N_{4}^{\varepsilon}},
$$
for $\varepsilon\in(0,a-1/2)$ small enough.
By using (\ref{GH-trilinear}), we have that the contribution of this case to the left-hand side of (\ref{i2}) has the desired estimate provided $a>0$ and $b>\frac12-\frac{\varepsilon}{3}$.

\underline{Case $(A_8)$.}
In this case region, we observe that $N_1\sim N_2,~N_3\sim N_4$ and $N_{14}\sim \max\{N_1,N_3\}$.
Then
\begin{eqnarray}\label{sa}
M(\overline{\tau},\overline{\xi})\lesssim \frac{\max\{K_1,K_2,K_3,K_4\}^{1/2}}{\max\{N_1,N_2,N_4,N_4\}^{2a-1/2}}.
\end{eqnarray}
Therefore, we have that the contribution of this case to the left-hand side of (\ref{i2}) has the desired estimate provided  $a>1/4$ and $b>1/2-(2a-1/2)/3$.

\underline{Case $(A_9)$.}
In this case region, we observe that
$$
N_1\sim N_2\sim N_3\sim N_4\gg N_{12}.
$$
The matter when $N_{14}\gtrsim N_1$ is reduced to the proof of Case $(A_8)$.
Indeed, in such a case, $M(\overline{\tau},\overline{\xi})$ satisfies (\ref{sa}), and we can repeat the argument presented above.
Hence it remains to consider the situation that (by symmetry)
$$
N_1\sim N_2\sim N_3\sim N_4\gg N_{12}\ge N_{14}.
$$
We shall consider this case in Lemma \ref{special case}.
As a sequel, we complete the proof of Lemma \ref{thm:tri-linear}.

\begin{lemma}\label{special case}
Let $s>4/9+a/9$ and $a>1/4$.
Given dyadic numbers\footnote{We use here a dyadic number to be a number $N=2^j$ where $j\in\mathbb{N}$.} $N$, suppose that for all $t\in\mathbb{R}$, $\mathrm{supp}{\cal F}u_j(t,\xi)\subset\{\xi\mid \langle \xi\rangle \sim N\}$ for $1\le j\le 3$.
There exist $c>0$ and $\varepsilon>0$ such that the following estimate holds
\begin{eqnarray}
& & \left(\int\sum_{\langle\xi\rangle\sim N}\frac{\langle\xi\rangle^{2a}}{\langle\tau+\xi^2\rangle}\left|\int_*\sum_*1_A(\xi,\xi_1,\xi_2,\xi_3){\cal F}u_1(\tau_1,\xi_1){\cal F}{\overline{\partial_xu_2}}(\tau_2,\xi_2){\cal F}{u_3}(\tau_3,\xi_3)\right|^2\right)^{\frac12}\label{sp1}\\
& \le & \frac{c}{N^{\varepsilon}}\left(\sum_{k=1}^3\|u_k\|_{X^{a,1/2}}\prod_{j=1,\ne k}^3\|u_j\|_{L^8_tH_x^s}
+\sum_{k=1}^3\|u_k\|_{L_t^{2}H_x^s}\prod_{j=1,\ne k}^3\left(\|u_j\|_{L_t^{\infty}H_x^s}+\|u_j\|_{L^8_tH_x^s}\right)\right),\nonumber
\end{eqnarray}
where
$$
A=\{(\xi,\xi_1,\xi_2,\xi_3)\in\mathbb{Z}^4\mid 0\ne \max\{|\xi_1+\xi_2|,|\xi_1-\xi|\}\ll N\}.
$$
\end{lemma}
\noindent
{\it Proof of Lemma \ref{special case}.}
It is convenient to use the notation from the proof of Lemma \ref{thm:tri-linear}.
We rewrite $\tau_4=-\tau$ and $\xi_4=-\xi$, and localize the frequencies $M_1\sim |\xi_{12}|,~M_2\sim |\xi_{14}|$ where $M_1$ and $M_2$ range over dyadic numbers.
Put $M_{min}=\min\{M_1,M_2\}$ and $M_{max}=\max\{M_1,M_2\}$.
Note that the identity (\ref{id}) implies that
$$
\max\{K_1,K_2,K_3,K_4\}\gtrsim |\xi_{12}||\xi_{14}|\sim M_1M_2.
$$
By symmetry, we analyze two cases
\begin{itemize}
\item[$(A_{91})$]
$K_4\gtrsim M_1M_2$,
\item[$(A_{92})$]
$K_1\gtrsim M_1M_2$.
\end{itemize}

\underline{Case $(A_{91})$.}
The convolution relation $\sum_{j=1}^4\xi_j=0$ implies $|\xi_{14}|=|\xi_{23}|$. 
By symmetry, we may assume $M_1\le M_2$.
Applying the Littlewoods-Paley inequality, it thus suffices to show that
\begin{eqnarray}
& & \sum_{M_1\le M_2\ll N}\frac{N^{a}}{M_1^{1/2}M_2^{1/2}}\left\|\sum_{*}\int_*1_{|\xi_{12}|\sim M_1}(\xi_1,\xi_2){\cal F}{u_1}(\tau_1,\xi_1){\cal F}{\overline{\partial_xu_2}}(\tau_2,\xi_2){\cal F}{u_3}(\tau_3,\xi_3)\right\|_{L_{\tau_4}^2\ell_{\xi_4}^2}\label{sp2}\\
& \lesssim & 
\frac{c}{N^{\varepsilon}}\|u_3\|_{L_t^{2}H_x^s}\prod_{j=1}^2\|u_j\|_{L_t^{\infty}H_x^s},\nonumber
\end{eqnarray}
where $M_j~(j=1,2)$ range over dyadic numbers with $M_1\le M_2\ll N$.
Undoing the Fourier transform with respect to time variable first and spatial variable next (we may assume $\widehat{u_j}(t,\xi)$ nonnegative for $t\in\mathbb{R}$), we bound the left-hand side of (\ref{sp2}) by
\begin{eqnarray}
& & \sum_{M_1\le M_2\ll N}\frac{N^{a}}{M_1^{1/2}M_2^{1/2}}\left\|{\cal F}_x^{-1}\left(\sum_*1_{|\xi_{12}|\sim M_1}(\xi_1,\xi_2)|\widehat{u_1}(t,\xi_1)||\widehat{\partial_x{u_2}}(t,\xi_2)|\right){\cal F}_x|\widehat{u_3}(t,\xi_3)|\right\|_{L_{t,x}^2}\nonumber\\
& \lesssim & \sum_{M_1\le M_2\ll N}\frac{N^{a}}{M_1^{1/2}M_2^{1/2}}\left\|{\cal F}_{x}^{-1}\left(\sum_*1_{|\xi_{12}|\sim M_1}(\xi_1,\xi_2)|\widehat{u_1}(t,\xi_1)||\widehat{\overline{\partial_xu_2}}(t,\xi_2)|\right)\right\|_{L_{t,x}^{\infty}}\|u_3\|_{L_{t,x}^2}.\label{sp3}
\end{eqnarray}
Using the Hausdorff-Young inequality, it follows that the first term in (\ref{sp3}) can be controlled by
\begin{eqnarray*}
& & 
\left\|{\cal F}_{x}^{-1}\left(\sum_*1_{|\xi_{12}|\sim M_1}(\xi_1,\xi_2)|\widehat{u_1}(t,\xi_1)||\widehat{\overline{\partial_x u_2}}(t,\xi_2)|\right)\right\|_{L_{t,x}^{\infty}} \\
& \lesssim & \left\|\sum_{\xi_1} 1_{|\xi|\sim M_1}(\xi_1,\xi-\xi_1)|\widehat{u_1}(t,\xi_1)||\widehat{\overline{\partial_x u_2}}(t,\xi-\xi_1)|\right\|_{L_t^{\infty}\ell_{\xi}^1}\\
& \lesssim & M_1 N\prod_{j=1}^2\|u_j\|_{L_t^{\infty}L_x^2}.
\end{eqnarray*}
Inserting this into (\ref{sp3}) and taking sum in $M_j$ reduced to showing that the left-hand side of (\ref{sp2}) is bounded by
$$
\frac{c}{N^{3s-a-1}}\log N\|u_3\|_{L_t^{2}H_x^s}\prod_{j=1}^2\|u_j\|_{L_t^{\infty}H^s}\le c\frac{c}{N^{\varepsilon}}\|u_3\|_{L_t^{2}H_x^s}\prod_{j=1}^2\|u_j\|_{L_t^{\infty}H^s},
$$
where $\varepsilon< 3s-a-1$.

\underline{Case $(A_{92})$.}
We use duality and will prove the following estimate
\begin{eqnarray}
& & 
\left|\int_{\Gamma_4}\sum_{\Gamma_4} {\cal F}{u_1}(\tau_1,\xi_1){\cal F}{\overline{\partial_x u_2}}(\tau_2,\xi_2){\cal F}{u_3}(\tau_3,\xi_3){\cal F}\overline{u_4}(\tau_4,\xi_4)\right|\label{sp5}\\
& \lesssim & N^{1+a-2s}M_{min}^{3/4}\|u_4\|_{X^{-a,b}}\|u_1\|_{L_{t,x}^2}\prod_{k=2}^3\|u_k\|_{L_t^8H_x^s},\nonumber
\end{eqnarray}
for $b>3/8$.
Notice that the above estimate implies that the contribution of this case to the left-hand side of (\ref{sp1}) is bounded by
\begin{eqnarray}\label{sp6}
cN^{1+a-2s}M_{min}^{3/4}\|u_1\|_{L_{t,x}^2}\prod_{k=2}^3\|u_k\|_{L_t^8H_x^s}.
\end{eqnarray}
Decompose each frequency with range $M_{min}$, namely
$$
\widehat{u_j}(t,\xi)=\sum_{k=1}^{N/M_{min}}\widehat{u_{j,k}}(t,\xi),
$$
where each $\widehat{u_{j,k}}$ has frequency with respect to $\xi$ within the range $M_{min}$.
In the region when $M_1\le M_2$ (the estimate in the case when $M_2>M_1$ is similar), there exists only one $l=l(k)$ of $u_{2,l}$ for each $k$ of $u_{1,k}$ (only one $m=m(n)$ of $u_{4,m}$ for each $n$ of $u_{3,n}$) such that
\begin{eqnarray*}
& & \sum_{\Gamma_4} \widehat{u_1}(t,\xi_1)\widehat{\overline{\partial_x u_2}}(t,\xi_2)\widehat{u_3}(t,\xi_3)\widehat{\overline{u_4}}(t,\xi_4)\\
& = & \sum_{k,n=1}^{N/M_1}\sum_{\Gamma_4}\widehat{u_{1,k}}(t,\xi_1)\widehat{\overline{\partial_x u_{2,l(k)}}}(t,\xi_2)\widehat{u_{3,n}}(t,\xi_3)\widehat{\overline{u_{4,m(n)}}}(t,\xi_4).
\end{eqnarray*}
Undoing the Fourier transform, it follows that the contribution of this case to the left-hand side is bounded by
\begin{eqnarray*}
& & c\left\|
\sum_{k,n=1}^{N/M_1}\|u_{1,k}(t)\|_{L_{x}^2}\|\partial_xu_{2,l(k)}(t)\|_{L_x^8}\|u_{3,n}(t)\|_{L_x^8}\|u_{4,m(n)}\|_{L_{x}^4}\right\|_{L_t^1}\nonumber\\
&\lesssim & \left\|\|u_{1,k}\|_{L_{x}^2}\right\|_{L_t^2\ell_k^2}
\left\|\|\partial_xu_{2,l(k)}(t)\|_{L_x^8}\right\|_{L_t^8\ell_k^2}
\left\|\|u_{3,n}(t)\|_{L_x^8}\right\|_{L_t^8\ell_n^2}
\left\|\|u_{4,m(n)}(t)\|_{L_x^4}\right\|_{L_t^4\ell_n^2},
\end{eqnarray*}
and by Minkowski's inequality, this is bounded by
\begin{eqnarray}\label{sp7}
c \|u_{1}\|_{L_{t,x}^2}
\left\|\|\partial_xu_{2,l(k)}(t)\|_{L_x^8}\right\|_{L_t^8\ell_k^2}
\left\|\|u_{3,n}(t)\|_{L_x^8}\right\|_{L_t^8\ell_n^2}
\left\|\|u_{4,m(n)}\|_{L_{t,x}^4}\right\|_{\ell_n^2}.
\end{eqnarray}
Since by Hausdorff-Young' and the decomposition of frequencies within the range $M_1$, we see that
$$
\|\partial_xu_{2,l(k)}(t)\|_{\ell_{k}^2L_x^8}
\lesssim N\|\widehat{u_{2,l(k)}}(t)\|_{\ell_k^2\ell_{\xi}^{8/7}}
\lesssim NM_1^{3/8}\|\widehat{u_{2,l(k)}}(t)\|_{\ell_{k,\xi}^2}
\lesssim NM_1^{3/8}\|u_2(t)\|_{L_x^2},
$$
and
$$
\|u_{3,n}(t)\|_{\ell_n^2L_x^8}\lesssim M_1^{3/8}\|u_3(t)\|_{L_x^2}.
$$
Moreover by (\ref{l6}), we see that
$$
\left\|\|u_{4,m(n)}\|_{L_{t,x}^4}\right\|_{\ell_n^2}\lesssim \left\|\|u_{4,m(n)}\|_{X^{0,1/2}}\right\|_{\ell_n^2}=\|u_4\|_{X^{0,b}},
$$
for $b>3/8$.
Inserting these estimates into (\ref{sp7}), we have that (\ref{sp7}) is bounded by 
\begin{eqnarray*}
c N^{1+a-2s}M_1^{3/4}\|u_1\|_{L^2_{t,x}}\|u_4\|_{X^{-a,b}}\prod_{j=2}^3\|u_j\|_{L_t^8H_x^s},
\end{eqnarray*}
which shows (\ref{sp5}).

Now since $\|u_1\|_{L^{2}_{t,x}}$ has two estimates, namely
\begin{eqnarray*}
\|u_1\|_{L_{t,x}^2}\lesssim\left\{
\begin{array}{l}
\frac{N^{-a}}{M_1^{1/2}M_2^{1/2}}\|u_1\|_{X^{a,1/2}},\\
N^{-s}\|u_1\|_{L_t^{2}H_x^s},
\end{array}
\right.
\end{eqnarray*}
Taking the sums in $M_2(\ge M_1)$ and $K_1\gtrsim M_1^{1/2}M_2^{1/2}$, the contribution of this case to the left-hand side of (\ref{sp1}) has two estimates
\begin{eqnarray}\label{interpolation-1}
cN^{1+a-3s}(\log  N)M_{min}^{3/4}\|u_1\|_{L_t^{2}H^s}\prod_{k=2}^3\|u_k\|_{L_t^{8}H_x^s},
\end{eqnarray}
and
\begin{eqnarray}\label{interpolation-2}
c\frac{N^{1-2s} \log N}{M_{min}^{1/4}}\|u_1\|_{X^{a,1/2}}\prod_{j=2}^3\|u_j\|_{L_t^8H_x^s}.
\end{eqnarray}
Interpolate (\ref{interpolation-1}) and (\ref{interpolation-2}) and taking the sum in $M_{min}\ll N$, it follows that the contribution of the case when $(A_{92})$ to the left-hand side of (\ref{sp1}) is bounded by
$$
\frac{c}{N^{\varepsilon}}\left(\|u_1\|_{L_t^{2}H_x^s}+\|u_1\|_{X^{a,1/2}}\right)\prod_{j=2}^3\left(\|u_j\|_{L_t^8H_x^s}+\|u_j\|_{L_t^{\infty}H_x^s}\right),
$$
for $s>4/9+a/9$, where $0<\varepsilon<(9s-a)/4-1$.
This completes the proof of Lemma \ref{special case}.
\qed

\noindent
{\it Proof of Lemma \ref{thm:tri-linear}.}
Now we return to the estimate for the case when $(A_9)$ in the proof of Lemma \ref{thm:tri-linear}.
Summing over dyadic number $N$ in Lemma \ref{special case}, we obtain that the contribution of the case $(A_9)$ to the left-hand side of (\ref{N-11-1}) is bounded by
$$
c\sum_{k=1}^3(\|u_k\|_{L_t^{2}H_x^s}+\|u_k\|_{X^{a,1/2}})\prod_{j=1,\ne k}^3(\|u_j\|_{L_t^8H_x^s}+\|u_j\|_{L_t^{\infty}H_x^s}),
$$
which leads to the result.
\qed

A proof similar to the one of Lemmas \ref{thm:tri-linear} and \ref{special case} allows us to prove the following lemma which is a variant of Lemma \ref{thm:tri-linear}.
\begin{lemma}\label{trilinear-2}
Let $4/9+a/9<s<1/2$ and $a>1/4$.
Then
\begin{eqnarray}
& & \left\|\langle\xi\rangle^a\frac{{\cal F}N_{11}(u_1,u_2,u_3)(\tau,\xi)}{\langle\tau+\xi^2\rangle}\right\|_{\ell_{\xi}^2L_{\tau}^1}\label{N-11-2}\\
& \lesssim &  \sum_{k=1}^3\|u_k\|_{X^{a,1/2}}\prod_{j=1,\ne k}^3\|u_j\|_{X^{a,b}}+ \sum_{k=1}^3\|u_k\|_{X^{a,1/2}}\prod_{j=1,\ne k}^3\|u_j\|_{L^8_tH_x^s}\nonumber\\
& & 
+\sum_{k=1}^3\|u_k\|_{L_t^{2}H_x^s}\prod_{j=1,\ne k}^3\left(\|u_j\|_{L_t^{\infty}H_x^s}+\|u_j\|_{L^8_tH_x^s}\right).\nonumber
\end{eqnarray}
\end{lemma}
\noindent
{\it Proof.}
We repeat the argument in the proof of Lemmas \ref{thm:tri-linear} and \ref{special case}.

In the region when $K_4\ll N_{12}N_{14}$, we show the required estimates from the proof of Lemmas \ref{thm:tri-linear} and \ref{special case} with subtle variation.
Indeed, using H\"older inequality in $\tau$, it follows that the contribution of this case to the left-hand side of (\ref{N-11-2}) is bounded by
\begin{eqnarray}\label{conv}
\|N_{11}(u_1,u_2,u_3)\|_{X^{a,-1/2+\varepsilon}}\left\|\frac{1}{\langle \tau+\xi^2\rangle^{1/2+\varepsilon}}\right\|_{\ell_{\xi}^{\infty}L_{\tau}^2}\lesssim\|N_{11}(u_1,u_2,u_3)\|_{X^{a,-1/2+\varepsilon}},
\end{eqnarray}
for $\varepsilon>0$.
By (\ref{id}) we see that
$$
\max\{K_1,K_2,K_3\}\gtrsim N_{12}N_{14}\gg K_4.
$$
In fact, using the trilinear estimate of (\ref{GH-trilinear}) with $b=\frac12-\varepsilon$, we have the following strong enough estimate
\begin{eqnarray}
& & \|N_{11}(u_1,u_2,u_3)\|_{X^{a,-1/2+\varepsilon}}\nonumber\\
& \lesssim & \sum_{k=1}^3\|u_k\|_{X^{a,1/2}}\prod_{j=1,\ne k}^3\|u_j\|_{X^{a,b}}+ \sum_{k=1}^3\|u_k\|_{X^{a,1/2}}\prod_{j=1,\ne k}^3\|u_j\|_{L^8_tH_x^s}\nonumber\\
& & 
+\sum_{k=1}^3\|u_k\|_{L_t^{2}H_x^s}\prod_{j=1,\ne k}^3\left(\|u_j\|_{L_t^{\infty}H_x^s}+\|u_j\|_{L^8_tH_x^s}\right),\label{t1}
\end{eqnarray}
which holds for $s>4/9+a/9,~a>(1+6\varepsilon)/4$ and some $3/8<b<1/2$.

On the other hand, in the region when $\max\{K_1,K_2,K_3\}\ll N_{12}N_{14}$, one notices that
\begin{eqnarray}\label{t2}
K_4\sim N_{14}N_{12}\lesssim \max\{N_1,N_2,N_3,N_4\}^2.
\end{eqnarray}
We review and change the proof of cases when $(A_j)$ for $j=2,3,4,5,7,8$ and $(A_{91})$.
For each case of $(A_j),~j=2,4,5,7,8$, by (\ref{t2}) we modify the bounds of $M(\overline{\tau},\overline{\xi})$ as follows
\begin{itemize}
\item[$(A_2)$]
$M(\overline{\tau},\overline{\xi})\lesssim N_1^{1/2-2a+\varepsilon}\max\{K_1,K_2,K_3,K_4\}$,
\item[$(A_4)$]
$M(\overline{\tau},\overline{\xi})\lesssim \frac{\max\{K_2,K_4\}^{1/2-\varepsilon}}{\min\{N_1,N_3\}^{2a-1/2}\max\{N_1,N_3\}^{1/2-\varepsilon}}$,
\item[$(A_5)$]
$M(\overline{\tau},\overline{\xi})\lesssim \frac{\max\{K_1,K_2,K_3,K_4\}^{1/2-\varepsilon}}{\max\{N_1,N_2,N_3,N_4\}^{\varepsilon}N_2^{2a-1/2-4\varepsilon}}$,
\item[$(A_7)$]
$M(\overline{\tau},\overline{\xi})\lesssim
\frac{\max\{K_1,K_2,K_3,K_4\}^{1/2-\varepsilon/2}}{ N_1^{a-\varepsilon/2}N_2^{a-\varepsilon}N_4^{\varepsilon/2}}$,
\item[$(A_8)$]
$M(\overline{\tau},\overline{\xi})\lesssim \frac{\max\{K_1,K_2,K_3,K_4\}^{1/2-\varepsilon}}{\max\{N_1,N_2,N_3,N_4\}^{2a-1/2-\varepsilon}}$,
\end{itemize}
for small $\varepsilon>0$.
It is not difficult to show that by the similar proof to the one in the case when $K_4\ll N_{12}N_{14}$ yields the result for $a>1/4$.

In the case when $(A_3)$, one notices that $K_4\sim N_2^2$.
We use the Littlewoods-Paley decompositions for $v_2$ as follows
$$
u_2(t,x)=\sum_{n}u_{2,n}(t,x),
$$
where $u_{2,n}(t,x)$ has the spatial Fourier support in the set $|\xi|\sim n$ for all $t\in\mathbb{R}$.
Also decompose $\langle \tau+\xi^2\rangle \sim m$ in the left-hand side of (\ref{N-11-2}).
By the restriction $K_4\sim N_2^2$, there exists one $m=m(n)$ for each $n$ of $u_{2,n}$ such that
\begin{eqnarray*}
& & \left\|\langle\xi\rangle^a\frac{{\cal F}N_{11}(u_1,u_2,u_3)(\tau,\xi)}{\langle\tau+\xi^2\rangle}\right\|_{\ell_{\xi}^2L_{\tau}^1}\\
& 
\lesssim &\sum_{n}\left\|1_{K_4\sim m(n)^2}(\tau,\xi)\langle\xi\rangle^s\frac{{\cal F}N_{11}(u_1,u_{2,n},u_3)(\tau,\xi)}{\langle\tau+\xi^2\rangle^{1/2}}\right\|_{L_{\tau}^2\ell_{\xi}^2}\left\|\frac{1_{K_4\sim m(n)^2}(\tau,\xi)}{\langle \tau+\xi^2\rangle^{1/2}}\right\|_{\ell_{\xi}^{\infty}L_{\tau}^2}.
\end{eqnarray*}
Notice that the second term in the right-hand side is bounded by a constant.
Reviewing the proof in the case when $(A_3)$ of the one of Lemma \ref{thm:tri-linear}, it suffices to prove that
\begin{eqnarray}
& & \sum_{n}\left|\int_{\Gamma_4}\sum_{\Gamma_4} 
M(\overline{\xi},\overline{\tau}){\cal F}{v_1}(\tau_1,\xi_1){\cal F}{\overline{v_{2,n}}}(\tau_2,\xi_2) {\cal F}{v_3}(\tau_3,\xi_3){\cal F}{\overline{v_{4,m(n)}}}(\tau_4,\xi_4)\right|\nonumber\\
&  \lesssim &\prod_{j=1}^4\|v_j\|_{X^{0,1/2}},\label{t3}
\end{eqnarray}
for any $v_4\in X^{0,1/2}$, where $v_j~(j=1,3)$ are same as in the proof of Lemma \ref{thm:tri-linear} and
$$
{\cal F}v_{2,n}(\tau,\xi)=\langle\xi\rangle^a{\cal F}u_{2,n}(\tau,\xi),~{\cal F}v_{4,m}(\tau,\xi)=1_{\langle\tau+\xi^2\rangle\sim m}(\tau,\xi)\langle\xi\rangle^{-a}{\cal F}v_4(\tau,\xi).
$$
Using the same proof in Lemma \ref{thm:tri-linear}, we have that the left-hand side of (\ref{t3}) is bounded by
$$
c\sum_{n}\|v_1\|_{X^{0,b}}\|v_{2,n}\|_{X^{0,b}}\|v_3\|_{X^{0,b}}\|v_{4,m(n)}\|_{X^{0,1/2}},
$$
for $b>1/2-a/3$, which is bounded by
$$
c\|v_4\|_{X^{0,1/2}}\prod_{j=1}^3\|v_j\|_{X^{0,b}},
$$
which yields the result for $a>0$. 

In the case when $(A_{91})$ in the proof Lemma \ref{special case}, we easily modify the proof as above and obtain that for small $\varepsilon>0$ the contribution of this case to the left-hand side of (\ref{N-11-2}) is bounded by
$$
c\sum_{M_1\le M_2\ll N}\frac{N^{1+a}M_1}{M_1^{1/2-\varepsilon}M_2^{1/2-\varepsilon}}\sum_{j=1}^3\|u_j\|_{L_t^2H^s}\prod_{k=1,\ne j}^3\|u\|_{L_t^{\infty}H^s},
$$
which is bounded by
$$
\frac{c}{N^{3s-a-1-3\varepsilon}}\sum_{j=1}^3\|u_j\|_{L_t^2H^s}\prod_{k=1,\ne j}^3\|u\|_{L_t^{\infty}H^s}.
$$
Then the result yields the desired estimate for $s>(a+1)/3$.
Therefore the proof is completed.
\qed

As a consequence of Lemmas \ref{thm:tri-linear} and \ref{trilinear-2}, we obtain the following proposition.
\begin{proposition}\label{prop-tri}
Let $4/9+a/9<s<1/2$ and $a>1/4$.
Then there exists $3/8<b<1/2$ such that
\begin{eqnarray*}
& & \|N_{11}(u_1,u_2,u_3)\|_{Z^a}\\
& \lesssim &  \sum_{k=1}^3\|u_k\|_{X^{a,1/2}}\prod_{j=1,\ne k}^3\|u_j\|_{X^{a,b}}+ \sum_{k=1}^3\|u_k\|_{X^{a,1/2}}\prod_{j=1,\ne k}^3\|u_j\|_{L^8_tH_x^s}\\
& & 
+\sum_{k=1}^3\|u_k\|_{L_t^{2}H_x^s}\prod_{j=1,\ne k}^3\left(\|u_j\|_{L_t^{\infty}H_x^s}+\|u_j\|_{L^8_tH_x^s}\right).\nonumber
\end{eqnarray*}
\end{proposition}

The next lemmas contains nonlinear estimates for $N_{12}[u],~N_{21}[u]$ and $N_{22}[u]$.

\begin{lemma}\label{lem-res}
For $s>0$, we have
\begin{eqnarray*}
\|N_{12}(u_1,u_2,u_3)\|_{Z^{s}}\lesssim \min_{1\le j\le 3}\|u_j\|_{L_t^{2}H^{\frac{s+1}{3}}}\prod_{k=1,\ne j}^3\|u_k\|_{L_t^{\infty}{\cal F}L^{\frac{s+1}{3},\infty}}.
\end{eqnarray*}
\end{lemma}
\noindent
{\it Proof.}
The proof is elementary, by using the fact that $\|f\|_{Z^s}\lesssim \|f\|_{L_t^{2}H^s},~L_x^2\hookrightarrow \ell_{\xi}^{\infty}$ and $|\xi|\langle\xi\rangle^s\lesssim \langle\xi\rangle^{\frac{s+1}{3}}$.
\qed

\begin{lemma}\label{lem-cubic}
For $s>0$ and $b>3/8$, we have
\begin{eqnarray*}
\|N_{21}[u]\|_{Z^{s}}\lesssim
\|u\|_{L_t^{2}H^s}\|u\|_{X^{0,b}}^2.
\end{eqnarray*}
\end{lemma}
\noindent
{\it Proof.}
By (\ref{l7}), we have that
$$
\|N_{21}[u]\|_{Z^{s}}\lesssim \|N_{21}[u]\|_{X^{s,-1/2+\varepsilon}}\lesssim\|\langle D_x\rangle^sN_{21}[u]\|_{L_{t,x}^{4/3}}.
$$
Using the Leibniz rule with fractional derivative and (\ref{l6}), it follows that 
$$
\|N_{21}[u]\|_{Z^{s}}\lesssim \|\langle D_x\rangle^su\|_{L_{t,x}^2}\|u\|_{L_{t,x}^4}^2,
$$
for small $\varepsilon>0$, which is bounded by $c\|u\|_{L_t^{2}H^s}\|u\|_{X^{0,b}}^2$, provided $b>3/8$.
\qed

\begin{lemma}\label{lem-quint}
For $s>1/4$ and $b>3/8$, we have
\begin{eqnarray}\label{tri-no}
\|N_{22}[u]\|_{Z^{s}}\lesssim\|u\|_{X^{s,b}}\|\langle D_x\rangle^{1/4}u\|_{L_t^8L_x^{8/3}}+\|u\|_{L_t^2H_x^s}\|u\|_{L_t^{\infty}H_x^{1/4}}^4.
\end{eqnarray}
\end{lemma}
\noindent
{\it Proof.}
We start by using (\ref{l7}),
$$
\|N_{22}[u]\|_{Z^{s}}\lesssim 
\|\langle D_x\rangle^sN_{22}[u]\|_{L_{t,x}^{4/3}}.
$$
By Leibniz rule with respect to fractional derivative, we have that the contribution of $|u|^4u$ term in $N_{22}[u]$ to (\ref{tri-no}) is estimated by
\begin{eqnarray}\label{tri-no-1}
c\|\langle D_x\rangle^su\|_{L_{t,x}^{4}}\|u\|_{L_{t,x}^8}^4.
\end{eqnarray}
Using (\ref{l4}) and Sobolev's inequality we conclude that (\ref{tri-no-1}) is bounded by
$$
c\|u\|_{X^{s,b}}\|\langle D_x\rangle^{1/4}u\|_{L_t^8L_x^{8/3}}^4.
$$

On the other hand, we have that the contribution of $\left(\int_0^{2\pi}|u(t,\theta)|^4\,d\theta\right)u$ to (\ref{tri-no}) is bounded by
$$
c\|u\|_{L_t^2H^s}\|u\|_{L_t^{\infty}L_x^4}^4.
$$
By Sobolev inequality, this is bounded by
$$
c\|u\|_{L_t^2H_x^s}\|u\|_{L_t^{\infty}H_x^{1/4}}^4.
$$
\qed

Finally, we shall attempt to localize the estimates in Proposition \ref{prop-tri}, Lemmas \ref{lem-res}, \ref{lem-cubic} and \ref{lem-quint}.
\begin{proposition}\label{ap}
Let $4/9+a/9<s<1/2$ and $a>1/4$.
Then there exists $\delta>0$ such that for any time $0<T<1$
\begin{eqnarray}\label{time-tri}
\|N_{11}[u]\|_{Z^a_T}\lesssim T^{\delta}\left(\|u\|_{Y^a_T}+\|u\|_{L_T^{\infty}H_x^s}\right)^3,
\end{eqnarray}
\begin{eqnarray}\label{time-11}
\|N_{12}[u]\|_{Z^a_T}\lesssim T^{\delta}\|u\|_{L_T^{\infty}H^s_x}^3,
\end{eqnarray}
\begin{eqnarray}\label{time-22}
\|N_{21}[u]\|_{Z^a_T}\lesssim T^{\delta}\|u\|_{Y_T^a}^3,
\end{eqnarray}
and
\begin{eqnarray}\label{time-33}
\|N_{22}[u]\|_{Z_T^a}\lesssim T^{\delta}\|u\|_{Y^a_T}^5.
\end{eqnarray}
\end{proposition}
\noindent
{\it Proof.}
For the sake of convenient, we only prove (\ref{time-tri}).
The estimates (\ref{time-11}), (\ref{time-22}) and (\ref{time-33}) follow using Lemmas \ref{lem-res}, \ref{lem-cubic} and \ref{lem-quint}, respectively.
Let $\widetilde{u}\in Y^a$ be such that $\widetilde{u}(t)=u(t)$ on $[-T,T]$.
We revisit the proof of Lemma \ref{thm:tri-linear} as well as  Proposition \ref{prop-tri}.
In cases when $(A_j),~1\le j\le 8$, we have that the contribution of these case to $\|N_{11}[u]\|_{Z^a_T}$ is bounded by
$$
\|N_{11}[\phi_T\widetilde{u}]\|_{Z^a}\lesssim \|\phi_T\widetilde{u}\|_{X^{a,1/2}}\|\phi_T\widetilde{u}\|_{X^{a,b}}^2,
$$
for some $3/8<b<1/2$.
By (\ref{l8}) and (\ref{l9}), this can be estimated as
\begin{eqnarray}\label{time-1}
\|N_{11}[\phi_T\widetilde{u}]\|_{Z^a}\lesssim T^{\delta}\|\widetilde{u}\|_{Y^a}^3.
\end{eqnarray}
In case when $(A_9)$, in a similar way as above, we have that the contribution of these case to $\|N_{11}[u]\|_{Z^a_T}$ is bounded by\footnote{In case when $(A_{92})$, we use the advantage in (\ref{sp5}) that by (\ref{l8}) $\|\phi_Tu_4\|_{X^{-a,b}}\lesssim T^{1/2-b-\varepsilon}\|u_4\|_{X^{-a,1/2}}$ for $3/8<b<1/2$ and $\varepsilon>0$.}
\begin{eqnarray}
& & cT^{b-3/8}\|\chi_T\widetilde{u}\|_{L_t^{\infty}H_x^s}^3+c\|\widetilde{u}\|_{X^{a,1/2}}\|\chi_T\widetilde{u}\|_{L_t^8H_x^s}^2+c\|\chi_T\widetilde{u}\|_{L_t^{\infty}H_x^s}\left(\|\chi_T\widetilde{u}\|_{L_t^{\infty}H_x^s}+\|\chi_T\widetilde{u}\|_{L_t^8H_x^s}\right)\nonumber\\
&\lesssim & T^{\delta}\left(\|\widetilde{u}\|_{Y^{a}}+\|u\|_{L_T^{\infty}H_x^s}\right)\|u\|_{L_T^{\infty}H_x^s}^2.
\label{time-2}
\end{eqnarray}
Therefore, by (\ref{time-1}) and (\ref{time-2}), we infer that
$$
\|N_{11}[u]\|_{Z_T^a}\lesssim T^{\delta}\left(\|\widetilde{u}\|_{Y^a}+\|u\|_{L_T^{\infty}H_x^s}\right)^3,
$$
which holds for any $\widetilde{u}$ satisfying $\widetilde{u}(t)=u(t)$ on $[-T,T]$.
Evaluate the infimum, then
$$
\|N_{11}[u]\|_{Z_T^a}\lesssim T^{\delta}\left(\|u\|_{Y^a_T}+\|u\|_{L_T^{\infty}H_x^s}\right)^3,
$$
which complets the proof of (\ref{time-tri}).
\qed

\section{Multilinear estimates II}\label{mul-1}
\indent

In this section we shall formulate and prove several preliminary estimates that are needed for the proof of Theorem \ref{thm-existence}.
\begin{lemma}[double mean value theorem]\label{dmvt}
Assume $f\in C^2(\mathbb{R})$ and that $\max\{|\eta|,|\lambda|\}\ll |\xi|$, then
$$
|f(\xi+\eta+\lambda)-f(\xi+\eta)-f(\xi+\lambda)+f(\xi)|\lesssim|f''(\theta)||\eta||\lambda|,
$$
where $|\theta|\sim |\xi|$.
\end{lemma}
\noindent
{\it Proof.}
See \cite[Lemma 2.3]{ckstt1}.

For $\overline{\xi}=(\xi_1,\xi_2,\xi_3,\xi_4)\in \mathbb{Z}^4\cap \Gamma_4$ with $\xi_{14}\xi_{34}\ne 0$, we let
\begin{eqnarray}\label{M_4}
M_4(\overline{\xi})=\frac{\xi_1\langle\xi_3\rangle^{2s}+\xi_2\langle\xi_4\rangle^{2s}+\xi_3\langle\xi_1\rangle^{2s}+\xi_4\langle\xi_2\rangle^{2s}}{\xi_{14}\xi_{34}}.
\end{eqnarray}

We have the following local estimate for $M_4$.

\begin{lemma}\label{lem:M_4-estimate}
Denote by $N_{(1)},~N_{(3)},~N_{(4)}$ the first, third, fourth biggest among $|\xi_j|~(1\le j\le 4)$, respectively.
Let $0<s<1/2$.
\begin{itemize}
\item[(i)]
If $N_{(1)}\lesssim \min\{|\xi_{14}|,|\xi_{34}|\}$ or $N_{(1)}\gg \max\{|\xi_{14}|,|\xi_{34}|\}$, then
\begin{eqnarray}\label{M4-1}
|M_4(\overline{\xi})|\lesssim
\langle N_{(1)}\rangle^{2s-1}.
\end{eqnarray}
\item[(ii)]
If $\max\{|\xi_{14}|,|\xi_{34}|\}\gtrsim N_{(1)}\gg \min\{|\xi_{14}|,|\xi_{34}|\}$, then
\begin{eqnarray}\label{M4-2}
|M_4(\overline{\xi})|\lesssim \langle N_{(3)}\rangle^{2s-1}.
\end{eqnarray}
\item[(iii)]
Assume that
$$
N_{(1)}=\max\{|\xi_1|,|\xi_3|\},~N_{(3)}=\max\{|\xi_2|,|\xi_4|\},~N_{(4)}=\min\{|\xi_2|,|\xi_4|\},
$$
or
$$
N_{(1)}=\max\{|\xi_2|,|\xi_4|\},~N_{(3)}=\max\{|\xi_1|,|\xi_3|\},~N_{(4)}=\min\{|\xi_1|,|\xi_3|\}.
$$
If $N_{(1)}\gg N_{(3)}$, then
\begin{eqnarray}\label{M4-3}
|M_4(\overline{\xi})|\lesssim
\langle N_{(3)}\rangle\langle N_{(1)}\rangle^{2s-2}.
\end{eqnarray}
\end{itemize}
\end{lemma}
\noindent
{\it Proof.}
Put $N_{(2)}$ the second biggest among $|\xi_j|$.
We have $N_{(1)}\sim N_{(2)}$ because of $\xi_1+\xi_2+\xi_3+\xi_4=0$ on $\Gamma_4$.
Since $\xi_{14}=-\xi_{23}$ and $\xi_{34}=-\xi_{12}$, by symmetry we may suppose $N_{(1)}=|\xi_1|$.

\underline{Case (i).}
We deal with the case when $|\xi_1|\lesssim \min\{|\xi_{14}|,|\xi_{34}|\}$ first.
In this case, we easily see that
$$
|M_4(\overline{\xi})|\lesssim \frac{\langle\xi_1\rangle^{2s+1}}{\langle\xi_1\rangle^2}\sim\langle\xi_1\rangle^{2s-1}.
$$
In the case when $|\xi_1|\gg \max\{|\xi_{14}|,|\xi_{34}|\}$, one has $|\xi_1|\sim|\xi_2|\sim|\xi_3|\sim|\xi_4|$ and $\xi_1\xi_4<0,~\xi_1\xi_2<0,~\xi_2\xi_3<0,~\xi_2\xi_4<0$.
We see that
\begin{eqnarray}
& & \langle\xi_1\rangle^{2s}\xi_3+\langle\xi_2\rangle^{2s}\xi_4+\langle\xi_3\rangle^{2s}\xi_1+\langle\xi_4\rangle^{2s}\xi_2\nonumber\\
& =  & \xi_{12}(\langle\xi_3\rangle^{2s}-\langle\xi_2\rangle^{2s})+\xi_{34}(\langle\xi_1\rangle^{2s}-\langle\xi_4\rangle^{2s})\nonumber\\
& & +\xi_{13}(\langle\xi_4\rangle^{2s}-\langle\xi_2\rangle^{2s})+\xi_{14}(\langle\xi_2\rangle^{2s}-\langle\xi_1\rangle^{2s})\label{1-1}\\
& & +\xi_1\langle\xi_1\rangle^{2s} +\xi_2\langle\xi_2\rangle^{2s} +\xi_3\langle\xi_3\rangle^{2s} +\xi_4\langle\xi_4\rangle^{2s}.\label{1-2}
\end{eqnarray}
By using mean value theorem, it follows that the term (\ref{1-1}) is bounded by
\begin{eqnarray*}
& & |\xi_{12}(\langle\xi_3\rangle^{2s}-\langle\xi_2\rangle^{2s})+\xi_{34}(\langle\xi_1\rangle^{2s}-\langle\xi_4\rangle^{2s})\\
& & +\xi_{13}(\langle\xi_4\rangle^{2s}-\langle\xi_2\rangle^{2s})+\xi_{14}(\langle\xi_2\rangle^{2s}-\langle\xi_1\rangle^{2s})|\\
& \lesssim & 
|\xi_{14}||\xi_{34}|\langle\xi_1\rangle^{2s-1}.
\end{eqnarray*}
Also using Lemma \ref{dmvt} (double mean value theorem), it follows that the term (\ref{1-2}) is bounded by
\begin{eqnarray*}
& & |\xi_1\langle\xi_1\rangle^{2s} +\xi_2\langle\xi_2\rangle^{2s} +\xi_3\langle\xi_3\rangle^{2s} +\xi_4\langle\xi_4\rangle^{2s}|\\
& = & 
|\xi_1\langle\xi_1\rangle^{2s}-(\xi_1-\xi_{14})\langle \xi_1-\xi_{14}\rangle^{2s}\\
& & -(\xi_1-\xi_{12})\langle\xi_1-\xi_{12}\rangle^{2s}+(\xi_1-\xi_{14}-\xi_{12})\langle\xi_1-\xi_{14}-\xi_{12} \rangle^{2s}|\\
& = & |\xi_{14}||\xi_{12}|\langle\xi_1\rangle^{2s-1}.
\end{eqnarray*}
These two estimates show that $|M_4(\overline{\xi})|\lesssim \langle\xi_1\rangle^{2s-1}$, which completes the proof in the case when (i).

\underline{Case (ii).}
In the case when $\max\{|\xi_{14}|,|\xi_{34}|\}\gtrsim |\xi_1|\gg \min\{|\xi_{14}|,|\xi_{34}|\}$, by symmetry, we may assume that $|\xi_{14}|\gtrsim|\xi_1|\gg|\xi_{34}|$ and $|\xi_1|\sim|\xi_2|\gtrsim\max\{|\xi_3|,|\xi_4|\}$.
Since
$$
M_4(\overline{\xi})=\frac{\xi_{12}\langle\xi_3\rangle^{2s}+\xi_{34}\langle\xi_1\rangle^{2s}+\xi_4(\langle\xi_2\rangle^{2s}-\langle\xi_1\rangle^{2s})}{\xi_{14}\xi_{34}}+\frac{\xi_2(\langle\xi_4\rangle^{2s}-\langle\xi_3\rangle^{2s})}{\xi_{14}\xi_{34}},
$$
and $\xi_{34}=-\xi_{12}$, by using mean value theorem, it follows that
$$
|M_4(\overline{\xi})|\lesssim\langle\xi_1\rangle^{2s-1}+\frac{|\langle\xi_4\rangle^{2s}-\langle\xi_3\rangle^{2s}|}{|\xi_{34}|}.
$$
For the second term in the right-hand side, we divide two cases that $|\xi_3|\sim|\xi_4|$ and that $\max\{|\xi_3|,|\xi_4|\}\gg\min\{|\xi_3|,|\xi_4|\}$.
If $|\xi_3|\sim|\xi_4|$, we again use the mean value theorem, while if $\max\{|\xi_3|,|\xi_4|\}\gg\min\{|\xi_3|,|\xi_4|\}$, we use $|\xi_{34}|\sim \max\{|\xi_3|,|\xi_4|\}$.
Then the second term is bounded by $c\min\{\langle \xi_3\rangle^{2s-1},\langle\xi_4\rangle^{2s-1}\}$, which completes the proof of the case (ii).

\underline{Case (iii).}
In this case, we may assume that
$$
N_{(1)}=|\xi_1|,~N_{(2)}=|\xi_3|,~N_{(3)}=|\xi_2|,~N_{(4)}=|\xi_4|,
$$
without loss of generality.
We rewrite
$$
M_4(\overline{\xi})=\frac{\xi_1(\langle\xi_1\rangle^{2s}-\langle\xi_3\rangle^{2s})}{\xi_{14}\xi_{34}}+\frac{\xi_{13}\langle\xi_1\rangle^{2s}}{\xi_{14}\xi_{34}}+\frac{\xi_2\langle\xi_4\rangle^{2s}+\xi_4\langle\xi_2\rangle^{2s}}{\xi_{14}\xi_{34}}.
$$
We apply the mean value theorem to the first term.
Since $|\xi_{34}|\sim |\xi_1|,~|\xi_{14}|\sim |\xi_1|,~|\xi_{13}|=|\xi_{24}|\lesssim |\xi_2|$, it follows that
$$
|M_4(\overline{\xi})|\lesssim|\xi_2|\langle\xi_1\rangle^{2s-2},
$$
which completes the proof in the case (iii).
\qed

We establish the following multilinear estimates.
\begin{lemma}\label{lem:m0}
For $1/4<s<1/2$, we have
\begin{eqnarray}
& & \left|\sum_{\Gamma_4}M_4(\overline{\xi})\left[\widehat{u}(t,\xi_1)\widehat{\overline{u}}(t,\xi_2)\widehat{u}(t,\xi_3)\widehat{\overline{u}}(t,\xi_4)\right]_{s=0}^{s=T}\right|\label{m-1}\\
& \lesssim &
\|u(T)\|_{H^{s/2}}^4+\|u_0\|_{H^{s/2}}^4.\nonumber
\end{eqnarray}
\end{lemma}
\noindent
{\it Proof.}
By Lemma \ref{lem:M_4-estimate}, we have
$$
|M_4(\overline{\xi})|\lesssim N_{(1)}^{s/2-1/8}N_{(2)}^{s/2-1/8}N_{(3)}^{s/2-3/8}N_{(4)}^{s/2-3/8}.
$$
Using Sobolev inequalities $L_x^8\hookrightarrow H_x^{3/8}$ and $L_x^{8/3}\hookrightarrow H_x^{1/8}$ along with the above inequality, we obtain the desired estimate.
\qed

\begin{lemma}\label{lem:m}
Let $s$ and $a$ with $1/4<a<s<\min\{1/2,3a/2\}$, and $0<T<1$.
Then there exists $\varepsilon>0$ such that
\begin{eqnarray}\label{m-2}
\left|\int_{0}^{T}\sum_{\Gamma_4}M_4(\overline{\xi})\widehat{u_1}(t,\xi_1)\widehat{\overline{u_2}}(t,\xi_2)\widehat{u_3}(t,\xi_3)\widehat{\overline{u_4}}(t,\xi_4)\,dt\right| \lesssim \|u_4\|_{Z^a}\prod_{j=1}^3\|u_j\|_{Y^{a}}.
\end{eqnarray}
\end{lemma}
\noindent
{\it Proof.}
By duality relation (\ref{duality-Y^s}), it suffices to show that
\begin{eqnarray}\label{m-3}
\left|\int_{-\infty}^{\infty}\sum_{\Gamma_4}M_4(\overline{\xi})\widehat{v_1}(t,\xi_1)\widehat{\overline{v_2}}(t,\xi_2)\widehat{v_3}(t,\xi_3)\widehat{\overline{v_4}}(t,\xi_4)\,dt\right|\lesssim \|v_4\|_{X^{a,-1/2}}\prod_{j=1}^3\|v_j\|_{Y^{a}},
\end{eqnarray}
and
\begin{eqnarray}\label{m-4}
\left|\int_{-\infty}^{\infty}\sum_{\Gamma_4}M_4(\overline{\xi})\left|\int_*{\cal F}{v_1}(\tau_1,\xi_1){\cal F}{\overline{v_2}}(\tau_2,\xi_2){\cal F}{v_3}(\tau_3,\xi_3)\right|\frac{w(-\xi_4)}{\langle\xi_4\rangle^a}\,d\tau_4\right|\lesssim \|w\|_{\ell^2_{\xi}}\prod_{j=1}^3\|v_j\|_{Y^{a}},
\end{eqnarray}
where $\tau_4=-(\tau_1+\tau_2+\tau_3)$.

First we consider (\ref{m-3}).
Use the dyadic partition $N_j\sim\langle\xi_j\rangle,~K_j\sim\langle\tau_j+(-1)^{j-1}\xi_j^2\rangle,~|\xi_{12}|\sim N_{12},~|\xi_{14}|\sim N_{14}$ as in the proof of Lemma \ref{thm:tri-linear}.
Since by (\ref{id}),
$$
\max\{K_1,K_2,K_3,N_{12}N_{14}\}\gtrsim K_4. 
$$
then separate the integral and sum of areas into following two cases
\begin{itemize}
\item[$(B_1)$]
$\max\{K_1,K_2,K_3\}\gtrsim K_4$,
\item[$(B_2)$]
$N_{12}N_{14}\gtrsim K_4$.
\end{itemize}

\underline{Case $(B_1)$.}
By symmetry and convolution strain, we may suppose $K_1\gtrsim K_4$ and $N_j\sim N_{(1)}$ for some $1\le j\le 3$.
If $N_1\sim N_{(1)}$, we have that from (i)-(ii) in Lemma \ref{lem:M_4-estimate},
$$
|M_4(\overline{\xi})|\lesssim N_1^aN_4^aN_2^{a-\frac12-\varepsilon}N_3^{a-\frac12-\varepsilon},
$$
for $s<3a/2$ and small $\varepsilon>0$.
Taking ${\cal F}^{-1}\frac{\langle\xi_4\rangle^a{\cal F}v_4}{\langle \tau_4-\xi_4^2\rangle^{1/2}}$ in $L_{t,x}^2$, it follows that this contribution to the left-hand side of (\ref{m-3}) is estimated by
\begin{eqnarray*}
 c\|v_4\|_{X^{a,-1/2}}\left\|\int_*\sum_*\langle\xi_1\rangle^a\langle\tau_1+\xi_1^2\rangle^{\frac12}|{\cal F}v_1(\tau_1,\xi_1)|\langle\xi_2\rangle^{a-\frac12-\varepsilon}|{\cal F}\overline{v_2}(\tau_2,\xi_2)|\langle\xi_3\rangle^{a-\frac12-\varepsilon}|{\cal F}v_3(\tau_3,\xi_3)|\right\|_{L_{\tau_4}^2\ell_{\xi_4}^2}.
\end{eqnarray*}
We use Plancherel's identity, Sobolev inequality, this is bounded by
$$
c\|v_4\|_{X^{a,-1/2}}\|v_1\|_{X^{a,1/2}}\prod_{j=2}^3\left\|{\cal F}^{-1}|{\cal F}v_j|\right\|_{L_t^{\infty}H^a}\lesssim \|v_4\|_{X^{a,-1/2}}\prod_{j=1}^3\|v_j\|_{Y^a},
$$
where we use Riemann-Lebesgue $\|{\cal F}_{\tau}^{-1}|{\cal F}_tv|\|_{L_t^{\infty}}\le \|{\cal F}_tv\|_{L_{\tau}^{1}}$.
On the other hand, if $N_1\not\sim N_{(1)}$, suppose $N_2\sim N_{(1)}$ and use
$$
|M_4(\overline{\xi})|\lesssim N_{2}^aN_{4}^aN_{1}^{a-\frac12-\varepsilon}N_{3}^{a-\frac12-\varepsilon}.
$$
Taking ${\cal F}^{-1}[\langle\xi_1\rangle^{a-\frac12-\varepsilon}\langle\tau_1+\xi_1^2\rangle^{\frac12}|{\cal F}v_1|]$ in $L_t^{2}L_x^{\infty}$, ${\cal F}[\langle\xi_2\rangle^a|\overline{v_2}|]$ in $L_t^{\infty}L_x^2$ and ${\cal F}^{-1}[\langle\xi_3\rangle^{a-\frac12-\varepsilon}|{\cal F}v_3|]$ in $L_{t,x}^{\infty}$, then we have that as same as above,  this contribution to the left-hand side of (\ref{m-3}) is estimated by
$$
c\|v_4\|_{X^{a,-1/2}}\prod_{J=1}^3\|v_j\|_{Y^a},
$$
for $s<3a/2$.

\underline{Case $(B_2)$.}
Notice that at least two of four $N_j$ are bigger than $cN_{(1)}$ for small constant $c>0$.
In the region when three of four $N_j$ are bigger than $cN_{(1)}$, by Lemma \ref{lem:M_4-estimate} (i)-(ii), we see that
$$
K_4^{1/2}|M_4(\overline{\xi})|\lesssim N_{(1)}^{2s}\lesssim \prod_{j=1}^4N_j^{a-\varepsilon},
$$
for $s<3a/2$ and small $\varepsilon>0$.
Using (\ref{GH-trilinear}), it follows that this contribution to the left-hand side of (\ref{m-3}) is bounded by
$$
c\|v_4\|_{X^{a,-1/2}}\left\|\prod_{j=1}^3{\cal F}^{-1}\langle\xi_j\rangle^{a-\varepsilon}{\cal F}v_j \right\|_{L_{t,x}^2}\lesssim \|v_4\|_{X^{a,-1/2}}\prod_{j=1}^3\|v_j\|_{X^{a,1/2}}.
$$
In other case when two of four $N_j$ are smaller than $cN_{(1)}$ for small constant $c>0$, separate the sum of area into two cases
\begin{itemize}
\item[$(B_{21})$]
$N_{(1)}\sim N_j$ and $N_{(2)}\sim N_k$ are occupied by a pairs of two odd or even numbers $j,k$, 
\item[$(B_{22})$]
otherwise.
\end{itemize}

In the subregion when $(B_{21})$, we see that $N_{12}N_{14}\lesssim N_{(1)}N_{(3)}$, which reduces that
$$
K_4^{1/2}|M_4(\overline{\xi})|\lesssim N_{(1)}^{1/2}N_{(3)}^{2s-1/2}\lesssim \prod_{j=1}^4N_j^{a-\varepsilon},
$$
for $s<3a/2$ and small $\varepsilon>0$.
In similar way to above, this contribution to the left-hand side of (\ref{m-3}) has the desired estimate.

In the subregion when $(B_{22})$, by Lemma \ref{lem:M_4-estimate} (iii), we see that
$$
K_4^{1/2}|M_4(\overline{\xi})|\lesssim N_{(1)}^{2s-1}N_{(3)}\lesssim \prod_{j=1}^4N_j^{a-\varepsilon},
$$
for $s<3a/2$ and small $\varepsilon>0$.
As above, this contribution to the left-hand side of (\ref{m-3}) has the desired estimate.

Let us prove the estimate (\ref{m-4}).
Writing
$$
a_{j}(\xi)=\langle\xi\rangle^a\int_{-\infty}^{\infty}|{\cal F}v_j(\tau,(-1)^{j-1}\xi)|\,d\tau,\quad 1\le j\le 3,
$$
and $a_{4}(\xi)=|w(\xi)|$, one can estimate the left-hand side of (\ref{m-4}) by
$$
c\sum_{\Gamma_4}\frac{|M_4(\overline{\xi})|}{\langle\xi_1\rangle^a\langle\xi_2\rangle^a\langle\xi_3\rangle^a\langle\xi_4\rangle^a}\prod_{j=1}^4a_{j}(\xi_j),
$$
which by $|M_4(\overline{\xi})|\lesssim 1,~a>1/4$ and Sobolev's inequality $L_x^4\hookrightarrow H^{a}_x$ with $a>1/4$, is bounded by
$$
c\prod_{j=1}^4\|a_j\|_{\ell_{\xi}^2}\lesssim \|w\|_{\ell_{\xi}^2}\prod_{j=1}^3\|v_j\|_{Y^a},
$$
as desired.
\qed

\section{A priori estimates}\label{sec:6}
\indent

In this section we prove the a priori estimates of solution that are needed for the proof of Theorem \ref{thm-existence}.

\subsection{$L_T^{\infty}H_x^s$ estimate}\label{H_x^s-estimates}
\indent

In this subsection we will derive a priori estimates in the $L_T^{\infty}H_x^s$ norm.

\begin{theorem}\label{thm:H^s-estimate}
Let $s$ and $a$ with $4/9+a/9<s<\min\{1/2,3a/2\}$ and $a>8/25$, and $u\in C^{\infty}(\mathbb{R},H^{\infty})\cap Y^s_{T}$ be global in time unique smooth solution to (\ref{dnls-ori})-(\ref{data}).
Then there exists a constant $\delta>0$ such that for all $N>1$, the following estimate holds for $|t|\le T$
\begin{eqnarray}
\|u(t)\|_{H^s} & \lesssim &  \|u_0\|_{H^s}(1+\|u_0\|_{H^s})+\|u_0\|_{L^2}\|u(t)\|_{H^s}\nonumber\\
& & +T^{\delta}\left(\|u\|_{L_T^{\infty}H_x^s}^3+\|u\|_{L_T^{\infty}H_x^s}^5+\|u\|_{Y^a_T}^3\right).\label{apr}
\end{eqnarray}
\end{theorem}
\noindent
{\it Proof.}
In order to discuss it, we first present a preliminary result.
Let $L[u](t)=\|u(t)\|_{H^s}^2$, so that
$$
L[u](t)=\sum_{\xi\in\mathbb{Z}}\langle\xi\rangle^{2s}|\widehat{u}(t,\xi)|^2=\sum_{\Gamma_2}m_2(\xi_1,\xi_2)\widehat{u}(t,\xi_1)\widehat{\overline{u}}(t,\xi_2),
$$
where $m_2(\xi_1,\xi_2)=(\langle\xi_1\rangle^{2s}+\langle\xi_2\rangle^{2s})/2$.
Note that
$$
\mathrm{Re}\sum_{\xi\in\mathbb{Z}}\langle\xi\rangle^{2s}\widehat{u}(t,\xi)\overline{\widehat{\partial_x^2u}(t,\xi)}=
c\sum_{\Gamma_2}(\xi_1^2-\xi_2^2)m_2(\xi_1,\xi_2)\widehat{u}(t,\xi_1)\widehat{\overline{u}}(t,\xi_2)=0,
$$
since $\xi_1^{2}-\xi_2^{2}$ vanishes on the hyperplane $\xi_1+\xi_2=0$.
Note also that
$$
\mathrm{Re}\sum_{\xi\in\mathbb{Z}}\langle\xi\rangle^{2s}\widehat{u}(t,\xi)\overline{\widehat{iN_{22}[u]}(t,\xi)}=\mathrm{Im}\sum_{\xi\in\mathbb{Z}} \frac{1}{2\pi}\xi\langle\xi\rangle^{2s}|\widehat{u}(t,\xi)|^4=0,
$$
and
$$
\mathrm{Re}\sum_{\xi\in\mathbb{Z}}\langle\xi\rangle^{2s}\widehat{u}(t,\xi)\overline{\widehat{iN_1[u]}(t,\xi)}=\mathrm{Re}\sum_{\xi\in\mathbb{Z}}\langle\xi\rangle^{2s}\widehat{u}(t,\xi)\overline{\widehat{iN_{11}[u]}(t,\xi)}.
$$
since real part is zero.
Thus
\begin{eqnarray}\label{dL}
\partial_tL[u](t)=\mathrm{Re}\sum_{\xi\in\mathbb{Z}}\langle\xi\rangle^{2s}\widehat{u}(t,\xi)\overline{\widehat{iN_{11}[u]}(t,\xi)}+\mathrm{Re}\sum_{\xi\in\mathbb{Z}}\langle\xi\rangle^{2s}\widehat{u}(t,\xi)\overline{\widehat{iN_2[u]}(t,\xi)}.
\end{eqnarray}
Now we can write the first term as
$$
\sum_{\xi\in\mathbb{Z}}\langle\xi\rangle^{2s}\widehat{u}(t,\xi)\overline{\widehat{iN_{11}[u]}(t,\xi)}=\frac{1}{2\pi i}\sum_{\scriptstyle \Gamma_4 \atop{\scriptstyle \xi_{12}\xi_{14}\ne 0}}\langle\xi_1\rangle^{2s}\xi_3\widehat{u}(t,\xi_1)\widehat{\overline{u}}(t,\xi_2)\widehat{u}(t,\xi_3)\widehat{\overline{u}}(t,\xi_4).
$$
Using symmetrization rules
\begin{itemize}
\item[(i)]
among two couples $\{\xi_1,\xi_3\}$ and $\{\xi_2,\xi_4\}$, namely $\{\xi_1,\xi_3\}=\{\xi_2,\xi_4\}$,
\item[(ii)]
between $\xi_1$ and $\xi_3$,
\item[(iii)]
between $\xi_2$ and $\xi_4$,
\end{itemize}
we compute
\begin{eqnarray}\label{modified-energy}
\mathrm{Re}\sum_{\xi\in\mathbb{Z}}\langle\xi\rangle^{2s}\widehat{u}(t,\xi)\overline{\widehat{iN_{11}[u]}(t,\xi)}=
c\sum_{\scriptstyle \Gamma_4 \atop{\scriptstyle \xi_{12}\xi_{14}\ne 0}}m_4(\overline{\xi})\widehat{u}(t,\xi_1)\widehat{\overline{u}}(t,\xi_2)\widehat{u}(t,\xi_3)\widehat{\overline{u}}(t,\xi_4),
\end{eqnarray}
for some constant $c$, where
$$
m_4(\overline{\xi})=\xi_1\langle\xi_3\rangle^{2s}+\xi_2\langle\xi_4\rangle^{2s}+\xi_3\langle\xi_1\rangle^{2s}+\xi_4\langle\xi_2\rangle^{2s}.
$$
Integrating with respect to $t$, we see that
\begin{eqnarray}
\|u(t)\|_{H^s}^2 & = &  \|u_0\|_{H^s}^2\nonumber\\
& & +c\int_0^t\sum_{\scriptstyle \Gamma_4 \atop{\scriptstyle \xi_{12}\xi_{14}\ne 0}}m_4(\overline{\xi})\widehat{u}(t',\xi_1)\widehat{\overline{u}}(t',\xi_2)\widehat{u}(t',\xi_3)\widehat{\overline{u}}(t',\xi_4)\,dt'\label{diff-1}\\
& & +c\int_0^t\mathrm{Re}\sum_{\xi\in\mathbb{Z}}\langle\xi\rangle^{2s}\widehat{u}(t',\xi)\overline{\widehat{iN_2[u]}(t',\xi)}\,dt'.\label{diff-2}
\end{eqnarray}
For (\ref{diff-1}), we will rewrite the ansatz $w=e^{-it\partial_x^2}u$, which implies that $i\partial_t u+\partial_x^2u=ie^{it\partial_x^2}\partial_t w$ and
$$
\widehat{u}(t,\xi_1)\widehat{\overline{u}}(t,\xi_2)\widehat{u}(t,\xi_3)\widehat{\overline{u}}(t,\xi_4)=e^{-2i\xi_{14}\xi_{34}t}\widehat{w}(t,\xi_1)\widehat{\overline{w}}(t,\xi_2)\widehat{w}(t,\xi_3)\widehat{\overline{w}}(t,\xi_4).
$$
Therefore using integration by parts it follows that
\begin{eqnarray*}
& & \int_0^t\sum_{\scriptstyle \Gamma_4 \atop{\scriptstyle \xi_{12}\xi_{14}\ne 0}}m_4(\overline{\xi})\widehat{u}(t',\xi_1)\widehat{\overline{u}}(t',\xi_2)\widehat{u}(t',\xi_3)\widehat{\overline{u}}(t',\xi_4)\,dt'\\
& = & -2i\sum_{\scriptstyle \Gamma_4 \atop{\scriptstyle \xi_{12}\xi_{14}\ne 0}}M_4(\overline{\xi})\left[e^{-2i\xi_{14}\xi_{34}t'}\widehat{w}(t',\xi_1)\widehat{\overline{w}}(t',\xi_2)\widehat{w}(t',\xi_3)\widehat{\overline{w}}(t',\xi_4)\right]_{t'=0}^{t'=t}\\
& & -2i\int_0^t\sum_{\scriptstyle \Gamma_4 \atop{\scriptstyle \xi_{12}\xi_{14}\ne 0}}M_4(\overline{\xi})e^{-2i\xi_{14}\xi_{34}t'}\partial_{t'}\left(\widehat{w}(t',\xi_1)\widehat{\overline{w}}(t',\xi_2)\widehat{w}(t',\xi_3)\widehat{\overline{w}}(t',\xi_4)\right)\,dt'\\
& = & F_1[u](t)+F_2[u](t).
\end{eqnarray*}
From $u=e^{it\partial_x^2}w$, we may thus
\begin{eqnarray*}
F_1[u](t)=c\sum_{\scriptstyle \Gamma_4 \atop{\scriptstyle \xi_{12}\xi_{14}\ne 0}}M_4(\overline{\xi})\left[\widehat{u}(t',\xi_1)\widehat{\overline{u}}(t',\xi_2)\widehat{u}(t',\xi_3)\widehat{\overline{u}}(t',\xi_4)\right]_{t'=0}^{t'=t}.
\end{eqnarray*}
and 
\begin{eqnarray*}
F_2[u](t)& = & c\sum_{k,l=1}^2\int_0^t\sum_{\scriptstyle \Gamma_4 \atop{\scriptstyle \xi_{12}\xi_{14}\ne 0}}\left[M_4(\overline{\xi})\prod_{(m,n)=(1,3)}\widehat{N_{kl}[u]}(t',\xi_m)\widehat{u}(t',\xi_n)\prod_{p=2,4}\widehat{\overline{u}}(t',\xi_p)\right]\,dt'\\
& & +c\sum_{k,l=1}^2\int_0^t\sum_{\scriptstyle \Gamma_4 \atop{\scriptstyle \xi_{12}\xi_{14}\ne 0}}\left[M_4(\overline{\xi})\prod_{(m,n)=(2,4)}\widehat{\overline{N_{kl}[u]}}(t',\xi_m)\widehat{\overline{u}}(t',\xi_n)\prod_{p=1,3}\widehat{u}(t',\xi_p)\right]\,dt'.
\end{eqnarray*}
Also for (\ref{diff-2}), we put
$$
F_3[u](t)=c\int_0^t\mathrm{Re}\sum_{\xi\in\mathbb{Z}}\langle\xi\rangle^{2s}\widehat{u}(t',\xi)\overline{\widehat{iN_2[u]}(t',\xi)}\,dt'.
$$
Then
$$
\|u(t)\|_{H^s}^2\le \|u_0\|_{H^s}^2+\sum_{j=1}^3F_j[u](t).
$$
Combining Proposition \ref{ap}, Lemma \ref{lem:m0}, Lemma \ref{lem:m}, (\ref{duality-Y^s}) and (\ref{t-inc}) with the arguments used in the proof of Proposition \ref{ap}, we see that there exists $\delta>0$ such that for $|t|\le T$
\begin{eqnarray}
\|u(t)\|_{H^s}^2& \lesssim &  \|u_0\|_{H^s}^2+\|u(t)\|_{H^{s/2}}^4+\|u_0\|_{H^{s/2}}^4\nonumber\\
& & +T^{\delta}\|u\|_{Y^a_T}^3\left(\|u\|_{L_T^{\infty}H_x^s}+\|u\|_{L_T^{\infty}H_s}^5+\|u\|_{Y_T^a}^3\right).\label{ap_1}
\end{eqnarray}
For the term $\|u(t)\|_{H^{s/2}}^4$, separating out spatial frequencies into high and low components and using $L^2$ conservation law, we have
\begin{eqnarray}\label{s/2}
\|u(t)\|_{H^{s/2}}\le \|P_{\le N}u(t)\|_{H^{s/2}}+\|P_{\ge N}u(t)\|_{H^{s/2}}\le N^{s/2}\|u_0\|_{L^2}+\frac{1}{N^{s/2}}\|P_{\ge N}u(t)\|_{H^s},
\end{eqnarray}
which by choosing $N^s=\langle \|u(t)\|_{H^s}/\|u_0\|_{L^2}\rangle$, is bounded by 
$$
\|u(t)\|_{H^{s/2}}\lesssim \|u(t)\|_{H^s}^{1/2}\|u_0\|_{L^2}^{1/2}.
$$
Finally, by inserting this into (\ref{ap_1}) and taking square root, the desired estimate (\ref{apr}) follows.
\qed

\subsection{$Y_T^{a}$ estimate}\label{sec:Y^a}
\noindent

Our result in this subsection is the following theorem.
\begin{theorem}
Let $s$ and $a$ with $4/9+a/9<s<\min\{1/2,3a/2\}$ and $a>8/25$, and 
$u\in C^{\infty}(\mathbb{R},H^{\infty})\cap Y^a_T$ be a time global solution to (\ref{dnls-ori})-(\ref{data}).
Then there exist constant $\delta>0$ and $\varepsilon>0$ with $\delta>\varepsilon$ such that the following estimate holds
\begin{eqnarray}\label{apri}
\|u\|_{Y_T^{a}}\lesssim T^{-\varepsilon}\|u_0\|_{H^a}+T^{\delta}\left(\|u\|_{L_T^{\infty}H^s}^3+\|u\|_{Y^a_T}^3+\|u\|_{Y^a_T}^5\right).
\end{eqnarray}
\end{theorem}
\noindent
{\it Proof.}
We consider the integral equation associated to (\ref{dnls-ori})-(\ref{data}).
Let $\widetilde{u}\in Y^{\infty}$ be such that $u(t)=\widetilde{u}(t)$ on $[-T,T]$.
Establishing the equation (\ref{dnls-ori})-(\ref{data}) in the Duhamel form, it follows that $u(t)$ and $\widetilde{u}(t)$ solve
$$
\phi_T(t)u(t)=\phi_T(t)e^{it\partial_x^2}u_0+\sum_{k,l=1}^2\phi_T(t)\int_0^te^{i(t-s)\partial_x^2}\phi_{T}(s)N_{kl}[\phi_T\widetilde{u}](s)\,ds,
$$
on $|t|\le T$.
Set the right-hand side by $\Phi(\widetilde{u})$,
$$
\Phi(\widetilde{u})(t)=\phi_T(t)e^{it\partial_x^2}u_0+\sum_{k,l=1}^2\phi_T(t)\int_0^te^{i(t-s)\partial_x^2}\phi_{T}(s)N_{kl}[\phi_T\widetilde{u}](s)\,ds.
$$
We observe that $\phi_T\widetilde{u}(t)=\chi_Tu(t)$ on $|t|\le T$ and
$$
\|u\|_{Y_T^a}\le \|\phi_Tu\|_{Y_T^a}\le \|\phi_Tu\|_{Y^a}.
$$
Using (\ref{l1}), (\ref{l2}), (\ref{l9}) and (\ref{l10}), it follows that
$$
\|u\|_{Y_T^a}\le \|\Phi(\widetilde{u})\|_{Y^a}\le cT^{-\delta}\|\phi\|_{H^a}+cT^{-\delta}\sum_{k,l=1}^2\|\phi_{T}N_{kl}[\phi_{T}\widetilde{u}]\|_{Z^{a}}.
$$
By Proposition \ref{ap}, we see that there exist positive constants $\delta>\varepsilon>0$ such that
$$
\|u\|_{Y^a_T}\lesssim 
T^{-\varepsilon}\|u_0\|_{H^a}+T^{2\delta}\left(\|\chi_Tu\|_{L_t^{\infty}H^s}^3+\|\phi_T\widetilde{u}\|_{Y^a}^3+\|\phi_T\widetilde{u}\|_{Y^a}^5\right).
$$
We use the fact that $\|\chi_T\widetilde{u}\|_{L_t^{\infty}H^s}=\|u\|_{L_T^{\infty}H^s}$, and take the infimum condition $\widetilde{u}=u$ on $|t|\le T$ to obtain
$$
\|u\|_{Y^a_T}\lesssim T^{-\varepsilon}\|u_0\|_{H^a}+T^{\delta}\left(\|u\|_{L_T^{\infty}H^s}^3+\|u\|_{Y^a_T}^3+\|u\|_{Y^a_T}^5\right),
$$
as desired.
\qed

\subsection{$\|P_{\ge N}u\|_{L_T^{\infty}H_x^s}$ estimate}\label{sec:PN}
\indent

Define smooth upside-down Fourier multiplier on the Fourier transform side as follows
$$
P_{\gtrsim N}u(\xi)={\cal }{\cal F}_x^{-1}\left[\psi_N\widehat{u}\right].
$$
Let us quickly review the proof of Theorem \ref{thm:H^s-estimate} in Section \ref{H_x^s-estimates}.
If we considered the a priori estimate of $\|P_{\gtrsim N}u(t)\|_{H^s}^2$, the multiplier $M_4(\overline{\xi})$ defined in (\ref{M_4}) would be replaced by
$$
\widetilde{M_4}(\overline{\xi})=\frac{\xi_1\langle\xi_3\rangle^{2s}\psi_N(\xi_3)^2+\xi_2\langle\xi_4\rangle^{2s}\psi_N(\xi_4)^2+\xi_3\langle\xi_1\rangle^{2s}\psi_N(\xi_1)^2+\xi_4\langle\xi_2\rangle^{2s}\psi_N(\xi_2)^2}{\xi_{14}\xi_{34}}.
$$
Then it is very convenient that one could use $\widetilde{M_4}(\overline{\xi})$ instead of $M_4(\overline{\xi})$ in the argument in Section \ref{mul-1}.
The following estimate would follow by a variant of the proceeding arguments:
\begin{eqnarray*}
\|P_{\gtrsim N}u(t)\|_{H^s}^2& \lesssim &  \|P_{\gtrsim N}u_0\|_{H^s}^2+\|P_{\gtrsim N/10}u(t)\|_{H^{s/2}}^2\|u(t)\|_{H^{s/2}}^2+\|P_{\gtrsim N/10}u_0\|_{H^{s/2}}^2\|u_0\|_{H^{s/2}}^2\\
& & +\frac{T^{\delta}}{N^{2\varepsilon}}\|u\|_{Y^a_T}^3\left(\|u\|_{L_T^{\infty}H_x^s}+\|u\|_{L_T^{\infty}H^s}^5+\|u\|_{Y_T^a}^3\right),
\end{eqnarray*}
for some $\varepsilon>0$.
Choosing $C>0$ large, and taking the square root (if needed reformulate $N$ with $N/10$), it follows that
\begin{eqnarray}
\|P_{\ge N}u(t)\|_{H^s} & \lesssim & \|P_{\ge CN}u_0\|_{H^s}\left(1+\|u_0\|_{H^s}\right)\nonumber\\
& & +\frac{1}{N^{s/2}}\|P_{\ge CN}u(t)\|_{H^{s}}\|u(t)\|_{H^{s/2}}\nonumber\\
& & +\frac{T^{\delta}}{N^{\varepsilon}}\left(
\|u\|_{Y^a_T}^3+\|u\|_{L_T^{\infty}H_x^s}^3+\|u\|_{L_T^{\infty}H^s}^5\right),\label{ap-2}
\end{eqnarray}
on $t\in [-T,T]$ and for all $N>1$.

\subsection{A priori estimates}
\indent

As a consequence of subsections \ref{H_x^s-estimates}, \ref{sec:Y^a} and \ref{sec:PN}, we shall show some a priori estimates for solutions of (\ref{dnls-ori})-(\ref{data}).

\begin{theorem}\label{thm:L-estimate}
Let $s$ and $a$ with $4/9+a/4<s<\min\{1/2,3a/2\}$ and $a>8/25$, and 
$u(t)\in C_t^{\infty}H_x^{\infty}$ be a unique time global solution to (\ref{dnls-ori})-(\ref{data}) with small $\|u_0\|_{L^2}$ norm.
Then there exist a positive time $T=T(\|u_0\|_{H^s})>0$ and positive constants $\varepsilon,~\eta$ such that
\begin{eqnarray}\label{apr-1}
\|u\|_{L_T^{\infty}H_x^s}+T^{\eta}\|u\|_{Y_T^{a}}\le C,
\end{eqnarray}
\begin{eqnarray}\label{apr-2}
\|P_{\ge N}u\|_{L_T^{\infty}H^s}\le C\|P_{\ge CN}u_0\|_{H^s}+\frac{C}{N^{\varepsilon}},
\end{eqnarray}
for all $N>1$, where constants $C$ depend only on $\|u_0\|_{H^s}$ and $T$.
\end{theorem}
\noindent
{\it Proof.}
From Theorem \ref{thm:H^s-estimate}, we have that there exists  $\varepsilon>0$ such that $\|u_0\|_{L^2}\le \varepsilon\ll 1$ and
$$
\|u\|_{L_T^{\infty}H^s_x}\lesssim \|u_0\|_{H^s}(1+\|u_0\|_{L^2})+T^{\delta/2}\left(\|u\|_{L_T^{\infty}H_x^s}+\|u\|_{L_T^{\infty}H_x^s}^5+\|u\|_{Y^a_T}^5\right),
$$
where the term $\|u(t)\|_{H^s}\|u_0\|_{L^2}$ in the right-hand side of (\ref{apr}) is absorbed by the term on the left-hand side of (\ref{apr}).
Combining this and (\ref{apri}), one can choose $T>0$ so small that the estimate (\ref{apr-1}) for some $C>0$, since by bootstrap and continuity arguments.

The estimate (\ref{apr-2}) follows by subsection \ref{sec:PN}.
\qed

\section{Proof of Theorem \ref{thm-existence}}\label{sec:energy}
\indent

We now prove Theorem \ref{thm-existence}.
Fix $M>0$ and $T>0$ to be chosen later.
We construct a solution by a compactness theorem.
Given $u_0\in H^s$, we choose $u_{0,n}\in  H^s$ satisfying $u_{0,n}\to u_0$ in $H^s$.
Let now $M>0$ so large with $\|u_{0,n}\|_{H^s}\le M$ and $\|u_0\|_{H^s}\le M$ for all $n$.
Let $u_n$ be the time global solution of (\ref{dnls-ori}) corresponding to the initial data $u_{0,n}$.
It follows from Theorem \ref{thm:L-estimate} that there exist $T'\in (0,T]$ and $C>0$ depending only on $\|u_0\|_{H^s}$ such that
$$
\|u_n\|_{L_{T'}^{\infty}H^s}+\|u_n\|_{Y_{T'}^a}\le C,
$$
and
$$
\|P_{\ge N}u_n\|_{L_{T'}^{\infty}H^s}\le C\|P_{\ge CN}u_{0,n}\|_{H^s}+\frac{C}{N^{\alpha}},
$$
for all $n\in\mathbb{N}$.
Passage to the limit and applying the compactness theorem, we deduce that there exists a solution $u$ of (\ref{dnls-ori})-(\ref{data}) satisfying
$$
u\in L^{\infty}([-T',T'];H^s)\cap Y_{T'}^a,
$$
$$
\lim_{n\to\infty}\left(\|u_n-u\|_{L_{T'}^{\infty}H^s}+\|u_n-u\|_{Y_{T'}^a}\right)=0,
$$
\begin{eqnarray}\label{ap1}
\|u\|_{L_{T'}^{\infty}H^s}+\|u\|_{Y_{T'}^a}\le C,
\end{eqnarray}
and
\begin{eqnarray}\label{ap2}
\|P_{\ge N}u\|_{L_{T'}^{\infty}H^s}\le C\|P_{\ge CN}u_{0}\|_{H^s}+\frac{C}{N^{\alpha}}.
\end{eqnarray}
Now we will prove $u\in C([-T',T'];\,H^s)$.
Let $N>0$ be so large.
We divide $u(t)$ into a low frequency group $P_{\le N}u(t)$ and a high frequency part $P_{\ge N}u(t)$.
Since by (\ref{ap1}), we have
$$
\|P_{\le N}u\|_{Y^s_{T'}}\le cN^{s-a}\|u\|_{Y^a_{T'}}\le CN^{s-a}.
$$
From $C([-T',T'];H^s)\hookrightarrow Y_{T'}^{s}$ in Remark \ref{embed}, it is easy see that $P_{\le N}u\in C([-T',T'];H^s)$.
Combining this with (\ref{ap2}), we obtain the estimate
\begin{eqnarray*}
\limsup_{t\to t_0}\|u(t)-u(t_0)\|_{H^s}& \le & \lim_{t\to t_0}\|P_{\le N}(u(t)-u(t_0))\|_{H^s}\\
& & +2\|P_{\ge N}u\|_{L_{T'}^{\infty}H^s}\\
& \lesssim & \|P_{\ge CN}u_{0}\|_{H^s}+\frac{1}{N^{\alpha}},
\end{eqnarray*}
for all $t_0\in[-T',T']$.
By letting $N\to \infty$, we conclude $\lim_{t\to t_0}\|u(t)-u(t_0)\|_{H^s}$.
Hence $u\in C([-T',T'];H^s)$, which completes the proof of Theorem \ref{thm-existence}.
\qed

\section{Proof of Theorem \ref{approx}}\label{sec:8}
\noindent

In this section we present the finite dimensional approximation of the solution to (\ref{dnls-ori}) in $H^s\cap {\cal F}L^{s_1,p}$ with $1/4<s<1/2<s_1$ and $2<p<4$.

Throughout this section, it is assumed that $\mu$ in (\ref{N_{21}}) is the function with respect to $t$, namely $\mu[u](t)=\|u(t)\|_{L_x^2}^2/2\pi$.

We first recall the finite dimensional approximation equation (\ref{fdnls})-(\ref{fdata}) that was derived in \cite{nors}.
A similar computation that in Section 3 would allow us to rewrite (\ref{fdnls})-(\ref{fdata}) as the following form
\begin{eqnarray}\label{fd-dnls}
i\partial_tu^N+\partial_x^2u^N=P_{\le N}
\sum_{k,l=1}^2N_{kl}[u^N],
\end{eqnarray}
\begin{eqnarray}\label{df-data}
u^N|_{t=0}=u_0^N=P_{\le N}u_0,
\end{eqnarray}
where we replace $\mu$ by $\mu_N=\|\phi_N\|_{L^2}^2/2\pi$ at the coefficient of the nonlinear term $N_{21}$ in (\ref{N_{21}}).

\subsection{Multilinear estimates III}
\noindent

Now we recall the trilinear estimates obtained by Gr\"unrock and Herr in \cite[Lemmas 2.4 and 2.5]{GH}.
\begin{lemma}[\cite{GH}]\label{GH-1}
Let $2\le r\le p\le 4$.
Then,
$$
\|N_{11}(u_1,u_2,u_3)\|_{{\cal Z}^{1/2,r}}\lesssim \|u_1\|_{{\cal X}^{1/2,1/2}_{p,2}}\|u_2\|_{{\cal X}^{1/2,1/2}_{r,2}}\|u_3\|_{{\cal X}^{1/2,1/2}_{p,2}}.
$$
\end{lemma}
We prove the following variant of Lemmas \ref{thm:tri-linear} and \ref{GH-1}. 
\begin{lemma}
Let $1/4<s<1/2$.
There exist $b<1/2,~2<p<4$ and $\varepsilon>0$ such that
\begin{eqnarray}
& & \|N_{11}(u_1,u_2,u_3)\|_{X^{s,-1/2}}  \lesssim  
\sum_{\{k,j,l\}=\{1,2,3\}}\|u_k\|_{X^{s,1/2}}\|u_j\|_{X^{1/4+\varepsilon,b}}\|u_l\|_{X^{1/4+\varepsilon,b}}
\label{2-tri}\\
& & +\sum_{\{k,j,l\}=\{1,2,3\}}\|u_k\|_{X^{1/4+\varepsilon,1/2}}\|u_j\|_{X^{s,b}}\|u_l\|_{X^{1/4+\varepsilon,b}} +\min_{1\le j\le 3}\|u_j\|_{X^{s,b}}\prod_{k=1,\ne j}^3\|u_k\|_{{\cal X}^{1/2,b}_{p,2}}.\nonumber
\end{eqnarray}
\end{lemma}
\noindent
{\it Proof.}
We repeat the proof of Lemma \ref{thm:tri-linear}.
Under the same notation as in the proof of Lemma \ref{thm:tri-linear}, we consider
\begin{eqnarray}\label{2-tri-1}
\left\|\frac{\langle\xi\rangle^s}{\langle\tau+\xi^2\rangle^{1/2}}
\sum_{\scriptstyle * \atop{\scriptstyle (\xi-\xi_1)(\xi-\xi_3)=0}}\int_*{\cal F}u_1(\tau_1,\xi_1)\xi_2{\cal F}\overline{u_2}(\tau_2,\xi_2){\cal F}u_3(\tau_3,\xi_3)\right\|_{L_{\tau}^2\ell_{\xi}^2},
\end{eqnarray}
and distinguish the integral and sum of the areas into nine cases $A_j,~1\le j\le 9$.

In the cases when $A_j,~1\le j\le 8$, the same proof as that in Lemma \ref{thm:tri-linear} shows that the contribution to these cases to (\ref{2-tri-1}) is bounded by
$$
c\sum_{\{k,j,l\}=\{1,2,3\}}\|u_k\|_{X^{s,1/2}}\|u_j\|_{X^{1/4+\varepsilon,b}}\|u_l\|_{X^{1/4+\varepsilon,b}}+c\sum_{\{k,j,l\}=\{1,2,3\}}\|u_k\|_{X^{1/4+\varepsilon,1/2}}\|u_j\|_{X^{s,b}}\|u_l\|_{X^{1/4+\varepsilon,b}},
$$
for $s>1/4$, which is the desired estimate.

On the other hand, in the case when $A_9$; since $\langle\xi\rangle^a|\xi_2|\sim \min_{1\le j\le 3}\langle\xi_j\rangle^a\prod_{k=1,\ne j}\langle\xi_k\rangle^{1/2}$ and $\langle\xi_j\rangle^{1/2}\sim \langle\xi_k\rangle^{1/2}$ for $1\le j,k\le 3$, we may freely rearrange the trilinear element functions $u_j,~1\le j\le 3$ in $N_1(u_1,u_2,u_3)$.
More precisely, it suffices to show that the contribution of this case to (\ref{2-tri-1}) is bounded by
$$
c\min_{1\le j\le 3}\|u_j\|_{X^{s,b}}\prod_{k=1,\ne j}^3\|u_k\|_{{\cal X}^{1/2,b}_{p,2}},
$$
but this follows from Lemma \ref{GH-1}\footnote{It was shown that \cite[Lemmas 2.4 and 2.5]{GH} hold with $b=1/2$. But in the case when $A_9$, the contribution of this case to (\ref{2-tri-1}) holds for the extremal case $b<1/2$.} by choosing $r=2$.
\qed

Analogous to Proposition \ref{ap}, we have shall need the following estimates on $N_{kl}[u]$ for $1\le j,k\le 2$.
\begin{proposition}\label{ap-1}
Let $1/4<s<1/2<s_1$ and $2<p<4$.
Then there exist $\delta,\varepsilon>0$ such that for any time $0<T<1$
\begin{eqnarray}\label{ap-11}
\|N_{11}[u]\|_{{\cal Z}^{s,p}_T}\lesssim T^{\delta}\|u\|_{Y^s_T}\|u\|_{Y^{1/4+\varepsilon}_T\cap{\cal Y}^{1/2,p}_T}^2,
\end{eqnarray}
\begin{eqnarray}\label{ap-12}
\|N_{12}[u]\|_{{\cal Z}^{s,p}_T}\lesssim T^{\delta}\|u\|_{{\cal Y}^{s,p}_T}\|u\|_{{\cal Y}^{1/2,p}_T}^2,
\end{eqnarray}
\begin{eqnarray}\label{ap-13}
\|N_{21}[u]\|_{{\cal Z}^{s,p}_T}\lesssim T^{\delta}\|u\|_{Y^s_T}\|u\|_{Y^0_T}^4,
\end{eqnarray}
and
\begin{eqnarray}\label{ap-14}
\|N_{22}[u]\|_{{\cal Z}^{s,p}_T}\lesssim T^{\delta}\|u\|_{Y^s_T}\|u\|_{Y^{1/4}_T}^4,
\end{eqnarray}
where $\|u\|_{Y^{1/4+\varepsilon}_T\cap{\cal Y}^{1/2,p}_T}=\|u\|_{Y^{1/4+\varepsilon}_T}+\|u\|_{{\cal Y}^{1/2,p}_T}$.
\end{proposition}
\noindent
{\it Proof.}
The estimates in (\ref{ap-13}) and (\ref{ap-14}) follow by Proposition \ref{ap}, where we use the fact that from $\ell^2\hookrightarrow \ell^p$, $Z^s\hookrightarrow {\cal Z}^{s,p}$ provided $p>2$.
In (\ref{ap-13}), we use $\mu[u](t)\lesssim \|u\|_{Y^0}^2$.
 
We will see that (\ref{2-tri}) implies (\ref{ap-11}) by the same argument as in the proof of Propositions \ref{prop-tri} and \ref{ap}.

Finally, for the estimate (\ref{ap-12}), by $L_t^2{\cal F}L^{s,p}\hookrightarrow {\cal Z}^{s,p}$ it follows that
$$
\|N_{12}[u]\|_{{\cal Z}^{s,p}}\lesssim \left\||\xi|\langle\xi\rangle^s|\widehat{u}(t,\xi)^2\widehat{u}(t,\xi)\right\|_{L_t^2\ell_{\xi}^p},
$$
which by $\ell^p\hookrightarrow \ell^{\infty}$, is bounded by
$$
c\|u\|_{{\cal Y}^{s,p}}\|u\|_{{\cal Y}^{1/2,p}}^2.
$$
Again the same argument as that in Proposition \ref{ap} gives a gain of $T^{\delta}$ factor, which completes the proof.
\qed

We easily have more general estimates as follows.
\begin{corollary}\label{ap-2}
Let $1/4<s<1/2<s_1$ and $2<p<4$.
Then there exists $\delta>0$ such that for any time $0<T<1$
\begin{eqnarray*}
\|N_{11}[u_1]-N_{11}[u_2]\|_{{\cal Z}^{s,p}_T}\lesssim T^{\delta}\|u_1-u_2\|_{Y^s_T}\max_{j=1,2}\|u_j\|_{Y^s_T\cap{\cal Y}^{1/2,p}_T}^2,
\end{eqnarray*}
\begin{eqnarray*}
\|N_{12}[u_1]-N_{12}[u_2]\|_{{\cal Z}^{s,p}_T}\lesssim T^{\delta}\|u_1-u_2\|_{{\cal Y}^{s,p}_T}\max_{j=1,2}\|u_j\|_{{\cal Y}^{1/2,p}_T}^2,
\end{eqnarray*}
\begin{eqnarray*}
\|N_{21}[u_1]-N_{21}[u_2]\|_{{\cal Z}^{s,p}_T}\lesssim T^{\delta}\|u_1-u_2\|_{Y^s_T}\max_{j=1,2}\|u_j\|_{Y^0_T}^4,
\end{eqnarray*}
and
\begin{eqnarray*}
\|N_{22}[u_1]-N_{22}[u_2]\|_{{\cal Z}^{s,p}_T}\lesssim T^{\delta}\|u_1-u_2\|_{Y^s_T}\max_{j=1,2}\|u_j\|_{Y^{1/4}_T}^4.
\end{eqnarray*}
\end{corollary}

\subsection{Proof of Theorem \ref{approx}}
\noindent

We now prove Theorem \ref{approx}.
Fix $1/4<s'<s<1/2<s_1'<s_1$ and $\varepsilon>0$ such that $s'>1/4+\varepsilon$ and $s_1'>1/2+\varepsilon$.

First we recall the local well-posedness result for (\ref{dnls-ori})-(\ref{data}) in $H^s\cap {\cal F}L^{s_1,p}$.
We define the set
$$
A=\left\{u \in Y^s_T\cap {\cal Y}_T^{s_1} \mid ~\|u\|_{Y_T^s\cap{\cal Y}_T^{s_1}}\le M\right\},
$$
equipped with the distance
$$
|||u_1-u_2|||_{s,s_1,T}=\|u_1-u_2\|_{Y^s_T}+\|u_1-u_2\|_{{\cal Y}^{s_1}_T},
$$
where $M>0$ and $T>0$ are chooses later.

For $u_0\in H^s\cap {\cal F}L^{s_1,p}$, we define the operator
$$
\Psi[u](t)=\phi(t)e^{it\partial_x^2}u_0+\sum_{k,l=1}^2\phi(t)\int_0^te^{i(t-t')\partial_x^2}\phi_T(s)N_{kl}[u](t')\,dt'.
$$
As a consequence of Lemma \ref{linear}, Lemma \ref{l-gain}, Proposition \ref{ap}, Proposition \ref{ap-1}, Corollary \ref{ap-2} combining with the argument in subsection \ref{sec:Y^a}, we deduce that there exist $\delta>0$ such that
$$
|||\Psi[u]|||_{s,s_1,T}\le c \|u_0\|_{H^s\cap {\cal F}L^{s_1,p}}+cT^{\delta}M^3(1+M^2),
$$
$$
|||\Psi[u_1]-\Psi[u_2]|||_{s,s_1,T}\le c T^{\delta}|||u_1-u_2|||_{s,s_1,T}M^2(1+M^2),
$$
for $u,~u_1,~u_2\in A$.
Setting $M=2c\|u_0\|_{H^s\cap{\cal F}L^{s_1,p}}$ and $T$ such that $2cT^{\delta}M^2(a+M^2)<1$, we have that $\Psi$ defines a contraction map on $A$.
Therefore the Cauchy problem (\ref{dnls-ori})-(\ref{data}) is well-posed in the time interval $[-T_0,T_0]$ where $T_0=\langle \|u_0\|_{H^{1/4+\varepsilon}\cap {\cal F}L^{1/2+\varepsilon,p}}\rangle^{-\varepsilon'}$ for some $\varepsilon'>0$.
Analogously we can prove that the Cauchy problem (\ref{fd-dnls})-(\ref{df-data}) is well-posed in the same time interval.

Next we observe that
\begin{eqnarray*}
\left(i\partial_t+\partial_x^2\right)\left[u(t)-u^N(t)\right] =  
\sum_{k,l=1}^2P_{>N}N_{kl}[u] +P_{\le N}\sum_{k,l=1}^2\left(N_{kl}[u]-N_{kl}[u^N]\right).
\end{eqnarray*}
Break the time interval $[0,T]$ into discrete intervals of size $T_0$, and put $t_j=jT_0,~1\le j\le T/T_0$.
From the estimates in Lemma \ref{linear}, Lemma \ref{l-gain}, Proposition \ref{ap}, Proposition \ref{ap-1}, Corollary \ref{ap-2}, it is easy to see that
\begin{eqnarray*}
|||u-u^N|||_{s',s_1',t_1}& \lesssim & \|P_{>N}u_0\|_{H^{s'}\cap {\cal F}L^{s_1',p}}+t_1^{\delta}|||P_{\gtrsim N}u|||_{s',s_1',t_1}M_1^2(1+M_1^2)\\
 & &+t_1^{\delta}|||u-u^N|||_{s',s_1',t_1}M_1^2(1+M_1^2),
\end{eqnarray*}
where $M_1=\max\{|||u|||_{s',s_1',t_1},|||u^N|||_{s,s_1,t_1}\}$.
From the local well-posedness theory as above, one has
$$
|||u|||_{s',s_1',t_1}\lesssim \|u_0\|_{H^{s'}\cap {\cal F}L^{s_1',p}}\le A,
$$
for $t_1=T_0\sim \langle\|u_0\|_{H^{s'}\cap {\cal F}L^{s_1',p}}\rangle^{-\varepsilon'}$.
Since $\|P_{>N}u_0\|_{H^{s'}\cap {\cal F}L^{s_1',p}}\le AN^{\max\{s'-s,s_1'-s_1\}}$, choosing $t_1>0$ small we obtain
$$
|||u-u^N|||_{s',s_1',t_1}\le cAN^{\max\{s'-s,s_1'-s_1\}}+\frac12|||u-u^N|||_{s',s_1',t_1},
$$
for some constant $c>0$, which yields
$$
|||u-u^N|||_{s',s_1',t_1}\le 2cAN^{\max\{s'-s,s_1'-s_1\}}.
$$
The iteration scheme can be used directly to obtain
\begin{eqnarray*}
\|u(t_j)-u^{N}(t_j)\|_{H^{s'}\cap {\cal F}L^{s_1',p}}& \lesssim &  |||u-u^N|||_{s',s_1',t_j}\\
 & \lesssim & 2^jAN^{\max\{s'-s,s_1'-s_1\}}\\
& \sim & 2^{\frac{t_j}{T_0}}AN^{\max\{s'-s,s_1'-s_1\}}\\
& \le & C_1\mbox{exp}[C_2(1+A)^{C_3}t_j]N^{\max\{s'-s,s_1'-s_1\}},
\end{eqnarray*}
as long as the right-hand side remains less than $1$.
This leads the result.
\qed

\end{document}